\documentclass{amsart}

\usepackage{amsmath}	
\usepackage{amsthm}
\usepackage{amssymb}
\usepackage{amsfonts}
\usepackage{graphicx}

\def\RR{\mathbb{R}}

\def\TT{\mathbb{T}}

\def\Int{\mbox{\rm{Int}}\;}

\def\Per {\mbox{\rm{Per}}}
\def\PA {\mbox{$\mathbb{P}$\rm{A}}}
\def\dist {\mbox{\rm{dist}}}
\def\irr {{\mbox{\rm{\scriptsize irr}}}}

\def\cO{\mathcal{O}}
\def\cF{\mathcal{F}}
\def\cS{\mathcal{S}}
\def\cN{{\mathcal N}}
\def\cG{\mathcal{G}}

\def\GA{{\mbox{GA}}}

\def\Im	 {\mbox{\rm{Im} }}
\def\diam{\mbox{\rm{diam}}}

\def\Vol{\mbox{\rm{Vol}}}

\def\-{{\setminus}}

\def\ra{{\rightarrow}}
\def\st{{\;|\;}}

\def\del{\partial}

\def\hsp{{\hspace{5mm}}}

\newcommand\cl[1]{\overline{#1}}

\newtheorem{theo}{Theorem}[section]
\newtheorem{prop}[theo]{Proposition}
\newtheorem{lemm}[theo]{Lemma}
\newtheorem{coro}[theo]{Corollary}

\title{Regular projectively Anosov flows on three-dimensional manifolds}
\author{Masayuki Asaoka}
\address{Department of Mathematics, Kyoto University\\
 Kitashirakawa Oiwakecho, Sakyo-ku\\
 606-8502 Kyoto, Japan}

\email{asaoka@math.kyoto-u.ac.jp}
\thanks{Partially supported by Grant-in-Aid for Encouragement of Young
 Scientists (B) and JSPS Postdoctral Fellowships for Research Abroad}

\keywords{projectively Anosov flows, bi-contact structures}
\subjclass{37D30, 57R30}

\begin{document}

\begin{abstract}
We give the complete classification of regular
 projectively Anosov flows on closed three-dimensional manifolds.
More precisely, we show that
 such a flow must be either an Anosov flow or decomposed into
 a finite union of $T^2 \times I$-models.
We also apply our method to rigidity problems of some group actions.
\end{abstract}
\maketitle

\section{Introduction}
%
%
\subsection{Regular projectively Anosov flows}
In \cite{Mi}, Mitsumatsu introduced {\it a bi-contact structure}
 on a three-dimensional manifold, {\it i.e.},
 a pair of mutually transverse positive and negative contact structures.
He observed that a three-dimensional Anosov flow naturally induces
 a bi-contact structure whose intersection
 as a pair of plane fields is tangent to the flow.
In general, the intersection of a bi-contact structure
 does not define an Anosov flow.
In fact, he showed that a bi-contact structure corresponds
 to {\it a projectively Anosov flow}, which is a generalization
 of an Anosov flow.
In \cite{ET},
 Eliashberg and Thurston also studied bi-contact structures
 and projectively Anosov flows
 ({\it conformally Anosov flows} in their book)
 from the viewpoint of confoliation theory.
They observed that a bi-contact structure
 naturally appears in a linear deformation of
 a foliation into contact structures.

A flow $\Phi=\{\Phi^t\}_{t \in \RR}$ on a three-dimensional manifold
 $M$ is called a {\it projectively Anosov} flow (or a $\PA$ flow)
 if it has no stationary points and admits
 a decomposition $TM=E^u +E^s$ by continuous plane fields
 such that
\begin{itemize}
 \item $E^u(z) \cap E^s(z)=T \Phi(z)$
 for any $z \in M$,
 where $T\Phi$ is the line field tangent to the orbits of $\Phi$,
 \item $D\Phi^t(E^\sigma(z))=E^\sigma(\Phi^t(z))$
 for any $\sigma \in \{u,s\}$, $z \in M$, and $t \in \RR$, and
 \item there exist two constants $C>0$ and $\lambda >1$
 such that
\begin{displaymath}
\label{eq:PA def}
 \|N\Phi^t|_{(E^s/T\Phi)(z)}\| \cdot
 \|(N\Phi^t|_{(E^u/T\Phi)(z)})^{-1}\| \leq C\lambda^{-t}
\end{displaymath}
 for any $z \in M$ and $t \geq 0$,
 where $N\Phi=\{N\Phi^t\}_{t \in \RR}$
 is the flow on $TM/T\Phi$ induced by $\Phi$.
\end{itemize}
We call the decomposition $TM=E^u+E^s$ {\it a $\PA$ splitting}.
If it satisfies stronger inequalities
\begin{displaymath}
\label{eq:Anosov def}
 \|N\Phi^t|_{(E^s/T\Phi)(z)}\|\leq C\lambda^{-t},\;
 \|(N\Phi^t|_{(E^u/T\Phi)(z)})^{-1}\| \leq C\lambda^{-t}
\end{displaymath}
 for any $z \in M$ and $t \geq 0$,
 then the flow is called an {\it Anosov flow}
 and the splitting is called a {\it weak Anosov splitting}
\footnote{It is different from but equivalent to
 the common definition of an Anosov flow
 as pointed out by Doering \cite[Proposition 1.1]{Do}.}.
We remark that a $\PA$ splitting
 is a {\it dominated splitting} on the whole manifold.
Such a splitting plays important roles in the modern theory of dynamical
 systems. See \cite{BDV} for example.

It is known that a $\PA$ splitting is always integrable.
However, the splitting is not smooth in general\footnote{A $\PA$ flow
 with a smooth $\PA$ splitting is called {\it regular}.
 However, we do not use the term `regular' in this sense
 since we use this term in other context below.}.
In fact, any orientable closed three-dimensional manifold
 admits a smooth $\PA$ flow,
 but no $\PA$ flow on the three-dimensional sphere
 admits a $C^1$ $\PA$ splitting.
See Theorems 4.2.6 and 4.3.1 in \cite{Mi2}.
From the viewpoint of confoliation theory,
 a $\PA$ flow with a smooth $\PA$ splitting
 corresponds to a linear deformation
 of a smooth foliation into contact structures
 whose derivative generates another smooth foliation
 (see \cite[Proposition 2.2.3]{ET}).

In \cite{Gh}, Ghys classified
 three-dimensional Anosov flows with a smooth weak Anosov splitting.
We say an Anosov flow $\Phi$ is {\it algebraic} if
 there exists a Lie group $G$,
 its cocompact lattice $\Gamma$,
 and a one-parameter subgroup $\{a^t\}_{t \in \RR}$ of $G$
 such that $\Phi$ is a flow on $\Gamma \backslash G$
 given by $\Phi^t(\Gamma g)=\Gamma(g \cdot a^t)$.
It is known that there are only two choices of $G$ up to covering:
\begin{enumerate}
\item The universal covering group $\widetilde{\mbox{PSL}}(2,\RR)$
 of the special linear group of $\RR^2$.
\item The semi-direct product $\RR \ltimes \RR^2$
 associated to a homomorphism $H:\RR \ra \mbox{GL}(2,\RR)$
 given by $H(t)(x,y)= (e^t x, e^{-t},y)$.
\end{enumerate}
In the former case, the algebraic Anosov flow can be identified
 with the geodesic flow on a closed surface with a hyperbolic metric
 up to finite cover.
In the latter case, the algebraic Anosov flow can be identified
 with the suspension flow of a hyperbolic toral automorphism.
\begin{theo}
[\cite{Gh}]
If an Anosov flow on a closed three-dimensional manifold
 admits a $C^2$ $\PA$ splitting,
 it is smoothly equivalent to an algebraic Anosov flow.
\end{theo}

It is natural to ask whether
 any $\PA$ flow with a smooth $\PA$ splitting
 is equivalent to an algebraic model or not.
In \cite{No}, Noda showed that
 if a $\PA$ flow on a $\TT^2$-bundle over $S^1$
 admits a smooth $\PA$ splitting
 and has an invariant torus, then it must be represented 
 as a finite union of so-called {\it $\TT^2 \times I$-models}.
Roughly speaking,
 a $\TT^2 \times I$-model is a flow on $\TT^2 \times [0,1]$
 which is transverse to $\TT^2 \times \{z\}$ for any $z \in (0,1)$
 and is equivalent to a linear flow on each boundary.
See \cite{No} for the precise definition.
In a series of papers, he and Tsuboi gave a classification
  for certain manifolds, which is summarized as follows.
\begin{theo}[\cite{No,No2,NT,Ts}]
If a $\PA$ flow on
 a Seifert manifold or a $\TT^2$-bundle over $S^1$
 admits a smooth $\PA$ splitting,
 then it is either an Anosov flow or
 represented as a finite union of $\TT^2\times I$-models.
\end{theo}
The author of this paper also approached the classification
 from another direction.
In \cite{As}, he showed that if a $\PA$ flow
 on {\it any} closed three-dimensional manifold
 admits a smooth $\PA$ splitting
 and all periodic orbits are hyperbolic,
 then it is equivalent to one of the above.

In \cite{No2}, Noda conjectured that
 the above is the complete list of three-dimensional
 $\PA$ flows with a smooth $\PA$ splitting.
The goal of this paper is an affirmative solution to the conjecture.
\begin{theo}
\label{thm:main theorem}
If a $\PA$ flow on a closed, connected, and three-dimensional manifold
 admits a $C^2$ $\PA$ splitting,
 then it is either an Anosov flow or
 represented as a finite union of $\TT^2\times I$ models.
\end{theo}
The theorem gives a solution to a conjecture
 posed by Mitsumatsu (Conjecture 4.3.3 in \cite{Mi2})
 immediately.
\begin{coro}
Any bi-contact structure associated to a $\PA$ flow
 with a smooth $\PA$ splitting
 consists of tight contact structures.
\end{coro}

We give the proof of Theorem \ref{thm:main theorem}
 in Sections \ref{sec:dichotomy} and \ref{sec:transitive}.
In Section \ref{sec:dichotomy}, we show a dichotomy
 on dynamics of a $\PA$ flow with a $C^2$ $\PA$ splitting.
Namely, either the flow is topologically transitive
 or the non-wandering set is the union of
 invariant tori with rotational dynamics.
It is not so hard to see that the latter implies that the flow
 is represented by $\TT^2 \times I$-models.
In Section \ref{sec:transitive},
 we show the former implies that the flow is Anosov.
It is done by proving the hyperbolicity of all periodic orbits.

\subsection{Foliations with a tangentially contracting flow}
Let $\cF$ be a codimension-one foliation
 on a three-dimensional manifold $M$.
We say a flow $\Phi$ is {\it tangentially contracting}
 with respect to $\cF$
 if there exist $C>0$ and $\lambda>1$
 such that $\|N\Phi^t|_{T\cF/T\Phi(z)}\| \leq C\lambda^{-t}$
 for any $z \in M$ and $t \geq 0$.
We apply the method developed in this paper
 to a classification of foliations
 which admit a tangentially contracting flow.
\begin{theo}
\label{thm:TC}
Let $M$ be a closed three-dimensional manifold
 and $\cF$ a $C^r$ codimension-one foliation on $M$
 with $r \geq 2$.
Suppose that $\cF$ admits a $C^r$ tangentially contracting flow $\Phi$.
Then, $\Phi$ is Anosov and $\cF$ is $C^r$-diffeomorphic to
 the weak stable foliation of an algebraic Anosov flow.
\end{theo}

We give two examples of group actions which
 induce a foliation with a tangentially contracting flow naturally.
The above theorem implies the rigidity of such actions.

\paragraph{Locally free actions of the affine group}
Let $\GA$ be the group
 of orientation preserving affine transformations of
 the real line $\RR$.
It is generated by two one-parameter subgroups
 $\{a^t\}_{t \in \RR}$ and $\{b^x\}_{x \in \RR}$
 with a relation $b^x \cdot a^t = a^t \cdot b^{\exp(-t)x}$.
We say an action $\rho:M \times  \GA \ra M$ on a manifold $M$
 is {\it locally free} if the isotropy subgroup
 $\{g \in \GA \st \rho(p,g)=p\}$ is discrete
 for any $p \in M$.
By $\cO(p,\rho)$, we denote the $\rho$-orbit
 $\{\rho(p,g) \st g \in \GA\}$ of $p \in M$.
If $\rho$ is of class $C^1$ and $M$ is closed,
 then the partition $\cO_\rho=\{\cO(p,\rho) \st p \in M\}$ is a foliation.
The flow $\{\rho(\cdot,a^t)\}_{t \in \RR}$ is
 tangentially contracting with respect to $\cO_\rho$.

In \cite{Gh2}, Ghys classified $C^r$ locally free action
 of $\GA$ on closed three-dimensional manifolds for $r \geq 2$
 up to $C^r$ conjugacy
 assuming the existence of an invariant volume.
Applying Theorem \ref{thm:TC} to $\cO_\rho$,
 we obtain a classification of {\it the orbit foliation} of actions
 {\it without} the assumption on an invariant volume.
\begin{theo}
Let $\rho$ be a $C^r$ locally free action of $\GA$ 
 on a closed three-dimensional manifold with $r \geq 2$.
Then, the orbit foliation of $\rho$ is $C^r$ diffeomorphic to
 the weak stable foliation of an algebraic Anosov flow.
\end{theo}
In the forthcoming paper \cite{As0},
 we will give a classification of the actions of $\GA$
 {\it up to smooth conjugacy}.

\paragraph{Actions of Fuchsian groups on the circle}
Let $\Gamma_g$ be the fundamental group of
 the oriented closed surface of genus $g \geq 2$.
We identify the circle $S^1$ with the real projective line.
It induces a projective structure to the circle.
We call an action $\Phi$ of $\Gamma_g$ on the circle {\it projective}
 if it preserves the projective structure.
We say two actions $\Phi_1$ and $\Phi_2$ of $\Gamma_g$ are
 $C^r$-conjugate if there exists a $C^r$ diffeomorphism $H$
 (or a homeomorphism if $r=0$) of $S^1$ such that
 $H(\Phi_1(\gamma,p))=\Phi_2(\gamma,H(p))$
 for any $\gamma \in \Gamma_g$ and $p \in S^1$.
In \cite{Gh}, Ghys proved the rigidity of projective actions.
\begin{theo}
[\cite{Gh}]
\label{thm:Fuchsian}
Let $\Phi:\Gamma_g \times S^1 \ra S^1$  be a $C^r$ action
 with $r \geq 3$.
Suppose that $\Phi$ is $C^0$-conjugate to a projective action.
Then, it is $C^r$-conjugate to a projective action.
\end{theo}
His proof can be divided into two steps.
The first step is to show that the suspension foliation
 of the action admits a tangentially contracting flow.
The second step is to construct a transverse projective
 structure of the foliation using the flow obtained in the first step.
The first step can be done even for $r=2$.
The second step also can be done even for $r=2$
 if the flow obtained in the first step is Anosov,
 as Ghys mentioned in Section 5 of \cite{Gh}.
Hence, Theorem \ref{thm:TC} implies the following improvement
 of the above theorem.
\begin{theo}
Theorem \ref{thm:Fuchsian} holds even for $r=2$.
\end{theo}
It is known that the theorem does not hold for $r=1$.
See \cite{Gh}.

\paragraph{Acknowledgments}
The author would like to thank Professor Takashi Inaba
 for letting the author know the stability theory 
 and the level theory of Cantwell and Conlon,
 which are necessary to prove Lemma \ref{lemma:semi-proper}.
This paper was partially written when the author stayed at
 Unit\'e de Math\'ematiques Pures et Appliqu\'ees,
 \'Ecole Normale Sup\'erieure de Lyon.
He is grateful to the members of UMPA, especially to Professor
 \'Etienne Ghys for their warm hospitality.
He would also like to thank an anonymous referee
 for many comments to improve the paper.

\section{A dichotomy on dynamics}
\label{sec:dichotomy}

%
%
\label{sec:section2}
In the rest of the article,
 we fix an orientable, closed, connected,
 and three-dimensional manifold $M$.
Let $\Phi$ be a $C^2$ $\PA$ flow on $M$
 with a continuous $\PA$ splitting $TM=E^u + E^s$.
For a compact $\Phi$-invariant set $\Lambda$, 
 we define {\it the stable set} $W^s(\Lambda)$
 and {\it the unstable set} $W^u(\Lambda)$ by
\begin{eqnarray*}
 W^s(\Lambda) &=&
 \left\{z \in M
 \st \lim_{t \ra + \infty}d(\Phi^t(z),\Lambda) = 0\right\}
\end{eqnarray*}
 and $W^u(\Lambda)=W^s(\Lambda;\Phi^{-1})$,
 where $\Phi^{-1}$ is the time-reverse of $\Phi$.

We call a $\Phi$-invariant torus $T$ is {\it normally attracting}
 if there exists $C>0$ and $\lambda>1$ such that
 $\|N\Phi^t|_{E^s/T\Phi(z)}\| \leq C\lambda^{-t}$
 for any $z \in T$ and $t \geq 0$.
By the existence of a $\PA$ splitting,
 our definition coincides with the usual definition.
It is known that if $T$ is a normally attracting invariant torus,
 then $W^s(T)$ is an open neighborhood of $T$ and is diffeomorphic
 to $\TT^2 \times \RR$.
We say an invariant torus $T'$ is {\it normally repelling}
 if $T'$ is normally attracting
 with respect to the time-reverse $\Phi^{-1}$.

Let $\Omega_*$ be the union of invariant embedded tori
 to which the restriction of $\Phi$ are topologically equivalent
 to a linear flow.
For $\rho \in \{u,s\}$,
 let $\Omega^\rho_*$ be the union of tori in $\Omega_*$
 tangent to $E^\rho$.
By the linearity of the flow on tori and
 the domination property of the splitting,
 we have $\Omega_*=\Omega^u_* \cup \Omega^s_*$,
 $\Omega^u_*$ is a union of normally attracting $\Phi$-invariant tori,
 and $\Omega^s_*$ is a union of normally repelling $\Phi$-invariant tori.
The exactly same argument as Proposition 3.9 of \cite{AR}
 shows that $\Omega_*$ consists of finite number of tori.

The aim of this section is to show the following dichotomy.
\begin{prop}
\label{prop:dichotomy}
If $\Phi$ admits a $C^2$ $\PA$ splitting, then either
\begin{enumerate}
\item $\Phi$ is topologically transitive, or
\item $M=W^s(\Omega^u_*) \cup \Omega^s_*
 =W^u(\Omega^s_*) \cup \Omega^u_*$.
\end{enumerate}
\end{prop}
The latter implies that $\Phi$ is equivalent to one of the known models.
\begin{prop}
\label{prop:T2*I}
In the latter case of Proposition \ref{prop:dichotomy},
 $\Phi$ is represented by a finite union of $\TT^2 \times I$-models.
\end{prop}
\begin{proof}
Fix  connected components $T_0$ of $\Omega^s_*$
  and $U$ of $W^u(T_0) \- T_0$.
Take an embedding $\psi:\TT^2 \times [0,1] \ra M$
 such that $\psi(\TT^2 \times 0)=T_0$, $\Im \psi \subset U \cup T_0$,
 and $\psi(\TT^2 \times 1)$ is transverse to the flow.
Put $T'=\psi(\TT^2 \times 1)$.
Since $W^s(\Omega^u_*)$ is a disjoint union
 of the stable sets of connected components of $\Omega^u_*$,
 we have $T' \subset W^s(T_1)$
 for some connected component $T_1$ of $\Omega^u_*$.

Let $U_1$ be the connected component of $W^s(T_1) \- T_1$
 which contains $T'$.
Since $T_1$ is normally attracting
 there exists an embedding
 $\psi_1:\TT^2 \times [0,1] \ra M$
 such that $\Im \psi_1 \subset T_1 \cup (U_1 \- \Im \psi)$,
 $\psi_1(\TT^2 \times 0)=T_1$,
 $\psi_1(\TT^2 \times 1)$ is transverse to the flow,
 and each $\Phi$-orbit contained in $U_1$
 intersects with $\psi_1(\TT^2 \times 1)$ exactly once.
Then, we can take a smooth positive function $\tau$ on $T'$
 such that $\Phi^{\tau(z)}(z) \in \psi_1(\TT^2 \times 1)$
 for any $z \in T'$.
It implies that 
\begin{equation*}
 \cl{U}=\cl{U_1}
 =\Im \psi \cup \Im \psi_1 \cup
 \{\Phi^t(z) \st z \in T',t \in [0,\tau(z)]\}
\end{equation*}
 is diffeomorphic to $\TT^2 \times [0,1]$
 and its boundary is $T_0 \cup T_1$.

Inductively, we obtain sequences $(T_n)_{n \geq 0}$
 and $(B_n)_{n \geq 0}$ of subsets of $M$ such that
 $T_n$ is a connected component of $\Omega_*$,
 $B_n$ is diffeomorphic to $\TT^2 \times [0,1]$,
 $\del B_n=T_n \cup T_{n+1}$, and $B_n \cap B_{n+1}=T_{n+1}$
 for any $n$.
Since $\Omega_*$ contains only finitely many tori,
 we have $T_n=T_m$ for some $n \neq m$.
It implies that $M$ is a $\TT^2$-bundle over $S^1$.
By Noda's classification \cite{No},
 $\Phi$ is represented by a finite union of $\TT^2 \times I$-models.
\end{proof}

The rest of this section is devoted to
 the proof of Proposition \ref{prop:dichotomy}.
In Subsection \ref{sec:hyp-like},
 we show the existence of the stable and unstable
 manifolds of invariant sets without irregular periodic points.
In Section \ref{sec:reduced dichotomy},
 we prove the dichotomy
 when all periodic points outside $\Omega_*$ are regular.
In the both of subsections,
 we assume only a weaker condition on
 the regularity of the $\PA$ splitting
 since the $C^2$ regularity of the splitting is too strong
 when we apply the results
 to a foliation with a tangentially contracting flow.
At last,
 in Subsection \ref{sec:local dynamics},
 we prove that all periodic points are regular
 for $\PA$ flow with a $C^2$ $\PA$ splitting.

%
%
\subsection{Hyperbolic-like behavior}
\label{sec:hyp-like}
Let $\Phi$ be a $C^2$ $\PA$ flow
 and $TM=E^u + E^s$ its $\PA$ splitting.
In this subsection,
 we do not assume the the splitting is of class $C^2$.

For $z \in M$,
 we define {\it the orbit} $\cO(z)$,
 {\it the $\alpha$-limit set} $\alpha(z)$,
 and {\it the $\omega$-limit set} $\omega(z)$ by
\begin{align*}
 \cO(z) & =\{\Phi^t(z) \st z \in \RR\},\\
 \alpha(z) & =\bigcap_{T>0}\cl{\{\Phi^t(z) \st t  \leq -T\}},\\
 \omega(z) & =\bigcap_{T>0}\cl{\{\Phi^t(z) \st t  \geq T\}}.
\end{align*}
We say a point $z \in M$ is {\it periodic} if
 there exists $T>0$ such that $\Phi^T(z)=z$.
The minimum of $\{t>0 \st \Phi^t(z)=z\}$
 is called {\it the period} of $z$.
We denote the set of periodic points of $\Phi$ by  $\Per(\Phi)$,
 and the non-wandering set of $\Phi$ by $\Omega(\Phi)$.

We say a periodic point $z_0$ is {\it $s$-regular}
 when there exists an embedded closed annulus $A$ tangent to $E^s$
 such that $\Phi^t(A) \subset \Int A$ for any $t>0$
 and $\bigcap_{t>0} \Phi^t(A)=\cO(z_0)$.
Similarly, we say a periodic point $z_0$ is {\it $u$-regular} when
 there exists an embedded closed annulus $A$ tangent to $E^u$
 such that $\Phi^{-t}(A) \subset \Int A$ for any $t>0$
 and $\bigcap_{t>0} \Phi^{-t}(A)=\cO(z_0)$.
We also say $z_0$ is {\it $\rho$-irregular} for $\rho \in \{u,s\}$
 if $z_0$ is not $\rho$-regular.
Let $\Per^\rho_{\irr}(\Phi)$ be
 the set of $\rho$-irregular periodic points.
Put $\Per_\irr(\Phi)=\Per^s_\irr(\Phi) \cup \Per^u_\irr(\Phi)$.
The aim of this subsection is
 to show the existence of the unstable manifolds
 for a compact invariant set
 which does not intersect with $\Omega_* \cup \cl{\Per_\irr(\Phi)}$.

Fix a continuous family $\{\phi_z\}_{z \in M}$ of $C^2$ embeddings
 of $[-1,1]^2$ into $M$ such that
 $\Im \phi_z$ is transverse to $T\Phi$
 and $\phi_z(0,0)=z$ for any $z \in M$.
We call $\{\phi_z\}_{z \in M}$ {\it a family of local cross sections}.
Let $r_z^t$ be the holonomy map of the orbit foliation of $\Phi$
 between $\Im \phi_z$ and $\Im \phi_{\Phi^t(z)}$
 along the path $\{\Phi^{t'}(z) \st t' \in [0,t]\}$.
We call $\{r_z^t\}_{(z,t) \in M \times \RR}$
 {\it the family of local returns} associated to $\{\phi_z\}_{z \in M}$.
For $\Delta>0$, put
 $D_\Delta(z)=\{z' \in \Im \phi_z \st d(z,z') \leq \Delta\}$,
 where $d(z,z')$ is the distance of $z,z' \in M$.
By the continuity of the family $\{\phi^t_z\}_{z \in M}$,
 there exists $\Delta_\phi>0$
 such that $r_z^t$ is well-defined on $D_{\Delta_\phi}(z)$
 for any $z \in M$ and $t \in [-1,1]$.

The splitting $TM/T\Phi=(E^s/T\Phi) \oplus (E^u/T\Phi)$
 defines projections $\pi^s$ and $\pi^u$ 
 from $TM$ to $E^s/T\Phi$ and $E^u/T\Phi$ respectively.
For $\alpha>0$,
 we say an embedded interval $I$ in $M$
 is {\it an $(E^s,\alpha)$-transversal} if
 $\|\pi^s(v)\| \leq \alpha\|\pi^u(v)\|$
 for any $z \in I$ and $v \in T_z I$.
Similarly,
 we say an embedded interval $I$ in $M$
 is {\it an $(E^u,\alpha)$-transversal} if
 $\|\pi^u(v)\| \leq \alpha\|\pi^s(v)\|$
 for any $z \in I$ and $v \in T_z I$.
For $\Delta>0$,
 an interval $I$ is called {\it a $(\Delta,E^s)$-interval}
 if it is an $(E^s,1)$-transversal
 and $r_z^t(I) \subset D_\Delta(\Phi^t(z))$ for any $t \geq 0$.
Similarly,
 an interval $I$ is called {\it a $(\Delta,E^u)$-interval}
 if it is an $(E^u,1)$-transversal
 and $r_z^{-t}(I) \subset D_\Delta(\Phi^{-t}(z))$ for any $t \geq 0$.
 
The next lemma is a variant of ``the Denjoy property'',
 which was proved by Arroyo and Rodrigues-Hertz in \cite{AR} 
 for flows without non-hyperbolic periodic points.
\begin{lemm}
\label{lemma:Denjoy}
There exists $\Delta_0>0$ such that
\begin{enumerate}
\item  the interior of any $(\Delta_0,E^s)$-interval
 contains a point $z$ such that
 $\omega(z)$ is a periodic orbit in $\Per_\irr^u(\Phi)$
 or is a torus in $\Omega_*^u$, and
\item  the interior of
 any $(\Delta_0,E^u)$-interval contains a point $z$ such that
 $\alpha(z)$ is a periodic orbit in $\Per_\irr^s(\Phi)$
 or is a torus in $\Omega_*^s$.
\end{enumerate}
\end{lemm}
\begin{proof}
We will show the former assertion
 since that the latter can be obtained in the same way.
In the proof of Proposition 4.2 in \cite{AR},
 Arroyo and Rodrigues-Hertz used
 the hyperbolicity of all periodic points only in the proof of
 Lemma 4.3 and 4.4 and the other part of the proof works
even if there are non-hyperbolic periodic points.
Hence, it is sufficient to see
 how to recover the proof of Lemmas 4.3 and 4.4 for our case.

Fix a $(\delta,E^s)$-interval $I$ which contains $z$.
Put $I_t=r_z^t(I)$ for $t \geq 0$.
Let $\{J_s\}_{s \geq 0}$ be the family of $(\Delta,E^s)$-intervals
 in the proof of Proposition 4.2 of \cite{AR},
 {\it i.e.}, the maximal one
 among families of $(\Delta,E^s)$-intervals satisfying
 $I_s \subset J_s$ and $r_z^{t-s}(J_s) \subset J_t$
 for any $t \geq s \geq 0$.
As shown in Lemma 4.1 of \cite{AR}
 (its proof does not require the hyperbolicity of periodic orbits),
 there is a uniform bound from below of the length of
 the local stable manifold of each point of $J_s$.
Let $J_s^\epsilon$ be the union of the local stable manifold
 of each point of $J_s$.
Lemma 4.3 of \cite{AR} deals with the case that
 $r_z^t(J_s^\epsilon) \subset J_s^\epsilon$
 for some $z \in M$ and $t>0$.
In this case, the exactly same argument as Lemma 4.3 of \cite{AR}
 shows that the forward orbit of some point in $\Int J_0$
 converges to a $u$-irregular periodic orbit.

Lemma 4.4 of \cite{AR} deals with the case that
 $\limsup|J_s|>0$ and
 $r_z^t(J_s^\epsilon) \cap \Per(\Phi) \neq \emptyset$
 for some $s$ and $t>0$.
In this case, we need to show the inclination lemma
 for an $(\Delta,E^s)$-interval $J_s$
 and a periodic point which is sufficiently close to $J_s$.
However, it is an easy consequence of
 the existence of the local stable manifolds
 with uniform length (Lemma 4.1 in \cite{AR}).
\end{proof}

Let $\Lambda$ be a compact $\Phi$-invariant set
 such that $\Lambda \cap (\Omega_* \cup \cl{\Per_\irr(\Phi)})=\emptyset$.
In the rest of the subsection,
 we will show that the stable and unstable manifolds
 are well-defined for any point of $\Lambda$.

For a subset $S$ of $M$ and $\delta>0$,
 we denote {\it the $\delta$-neighborhood}
 $\{p \in M \st \inf_{q \in S} d(p,q) \leq \delta\}$
 by $\cN_\delta(S)$.
Fix $0<\Delta_1<\Delta_0$ such that
 $\cN_{\Delta_1}(\Lambda) \cap
  \left(\Omega_* \cup \cl{\Per_\irr(\Phi)}\right) = \emptyset$.
By the center-unstable manifold theorem,
 there exist constants $0<\delta_1<\delta_2<\Delta_1$
 and a continuous family $\{W^{cu}_{loc}(z)\}_{z \in M}$
 of $C^2$ $(E^s,1)$-transversals such that
\begin{itemize}
\item $z \in W^{cu}_{loc}(z) \subset D_{\delta_2}(z)$
 and $\del W^{cu}_{loc}(z) \subset \del D_{\delta_2}(z)$
 for any $z \in M$,
\item $W^{cu}_\delta(z)=W^{cu}_{loc}(z) \cap D_\delta(z)$ is an interval
 for any $0 <\delta<\delta_2$ and $z \in M$, and
\item $r_z^{-t}(W^{cu}_{\delta_1}(z)) \subset W^{cu}_{\delta_2}(\Phi^t(z))$
 for any $0 \leq t \leq 1$ and $z \in M$.
\end{itemize}

\begin{prop}
\label{prop:stable manifolds 1}
For any given $\delta>0$,
 there exists $\epsilon_1>0$ such that
\begin{enumerate}
 \item $r_z^{-t}(W^{cu}_{\epsilon_1}(z)) \subset W^{cu}_{\delta}(\Phi^t(z))$
 for any $z \in \Lambda$ and $t \geq 0$, and
 \item $\lim_{t \ra \infty}
 \left(\sup_{z \in M}|r_z^{-t}(W^{cu}_{\epsilon_1}(z))|\right)=0$.
\end{enumerate}
\end{prop}
\begin{proof}
Without loss of generality, we may assume $\delta<\delta_1$.
If the lemma does not hold, then there exist sequences
 $\{\epsilon_k>0\}_{k \geq 1}$, $\{t_k >0\}_{k \geq 1}$,
 and $(z_k \in \Lambda)_{k=0}^\infty$ such that
\begin{itemize}
\item $\lim_{k \ra \infty}\epsilon_k=0$,
\item $r_{z_k}^{-t}(W^{cu}_{\epsilon_k}(z_k))
 \subset W^{cu}_{\delta}(\Phi^{-t}(z_k))$
 for any $k \geq 0$ and $0 \leq t \leq t_k$, and
\item $\limsup_{k \ra \infty}|r_{z_k}^{-t_k}(W^{cu}_{\epsilon_k}(z_k))|>0$.
\end{itemize}
By the first and the last items,
 we have $\lim_{k \ra \infty} t_k=\infty$.
By taking subsequences, we may assume that
 $\Phi^{-t_k}(z_k)$ converges to a point $z_*$ of $\Lambda$
 and $r_{z_k}^{-t_k}(W^{cu}_{\epsilon_k}(z_k))$
 converges to an interval $I \subset W^{cu}_\delta(z_*)$
 with positive length. 
Then, $I$ is a $(\delta,E^s)$-interval.
Since $z_*$ is a point of  $I \cap \Lambda$,
 we have $\omega(z) \subset \cN_\delta(\Lambda)$
 for any $z \in I$,
 and hence,
 $\omega(z) \cap \left(\Omega_*(\Phi) \cup \cl{\Per_\irr(\Phi)}\right)
 =\emptyset$.
However, it contradicts Lemma \ref{lemma:Denjoy}.
\end{proof}
\begin{coro}
\label{cor:W-cu tangent}
For any sufficiently small $\epsilon>0$,
 $W^{cu}_{\epsilon}(z)$ is tangent to $E^u$ for any $z \in \Lambda$.
\end{coro}
\begin{proof}
Fix $\delta>0$ and
 take $\epsilon_1>0$ in Proposition \ref{prop:stable manifolds 1}.
For $z \in \Lambda$, $z' \in W^{cu}_{\epsilon_1}(z)$, and $t \geq 0$,
 let $\alpha(z,z',t)$ be the angle
 between $T_{r_z^{-t}(z')}r_z^{-t}(W^{cu}_{\epsilon_1}(z))$
 and $E^u(r_z^{-t}(z'))$.
By the domination property of the $\PA$ splitting,
 there exist $C>0$ and $\lambda>1$ such that
 $\alpha(z,z',t) \geq C\lambda^t\alpha(z,z',0)$.
On the other hand, the continuity of the family
 $\{W^{cu}_{\epsilon_1}(z)\}_{z \in M}$ and the proposition
 implies that $\alpha$ is bounded as a function of $z,z'$ and $t$.
Hence, we have $\alpha(z,z',0)=0$ for any $z$ and $z'$.
\end{proof}
\begin{lemm}
\label{lemma:E-u interval} 
There exists $\Delta_2 \in (0,\Delta_1/2)$ which satisfies
 the following property:
 If an $(E^u,1)$-transversal $I$
 is contained in $D_{\Delta_2}(z)$ for some $z \in M$
 and satisfies $r_z^{-t}(\del I) \subset D_{\Delta_2}(\Phi^{-t}(z))$
 for any $t \geq 0$,
 then it is a $(\Delta_1/2,E^u)$-interval.
\end{lemm}
\begin{proof}
Since $TM=E^s +E^u$ is a $\PA$ splitting,
 there exists $\alpha>0$ such that
 if $I$ is an $(E^u,1)$-transversal
 and $r_z^{-t}(I)$ is well-defined for some $z \in M$ and $t \geq 0$,
 then $r_z^{-t}(I)$ is an $(E^u,\alpha)$-transversal.
By the uniform transversality of
 $(E^u,\alpha)$-transversals to $E^u$,
 we can take $\delta \in (0,\Delta_1/2)$ and $\beta>1$ such that
 $|J| \leq \beta \cdot \diam(\del I)$
 for any $z \in M$ and any $(E^u,\alpha)$-transversal
 with $J \subset D_{\delta}(z)$. 
Put $\Delta_2=\delta/4\beta$.
We remark that $\Delta_2<\delta/4<\Delta_1/8$.

Let $I$ be an $(E^u,1)$-transversal
 contained in $D_{\Delta_2}(z)$ for some $z \in M$ such that
 $r_z^{-t}(\del I) \subset D_{\Delta_2}(\Phi^{-t}(z))$ for any $t \geq 0$.
It is sufficient to show that
\begin{equation*}
t_0=\sup\{t_1 \geq 0 \st r_z^{-t}(I) \subset D_{\delta}(\Phi^{-t}(z)),
 |r_z^{-t}(I)| \leq \delta \mbox{ for any } t \in [0,t_1]\}
\end{equation*}
 is infinite.
Suppose that $t_0$ is a finite number.
Since $r_z^{-t}(I)$ is an $(E^u,\alpha)$-transversal
 and $r_z^{-t}(\del I) \subset D_{\Delta_2}(\Phi^{-t}(z))$
 for $0 \leq t \leq t_0$, we have
\begin{equation}
\label{eqn:E-u interval 1} 
|r_z^{-{t_0}}(I)|
  \leq \beta \cdot \diam (r_z^{-{t_0}}(\del I))
 \leq 2 \beta \Delta_2 = \delta/2.
\end{equation}
It implies $|r_z^t(I)| <\delta$
 for any $t$ sufficiently close to $t_0$.
By the inclusion $r_z^{-t_0}(\del I) \subset D_{\Delta_2}(\Phi^{-t_0}(z))$
 again, the inequality (\ref{eqn:E-u interval 1}) implies
\begin{equation}
\label{eqn:E-u interval 2} 
r_z^{-t_0} (I) \subset D_{(\delta/2)+\Delta_2}(\Phi^{-{t_0}}(z))
 \subset D_{(3/4)\delta}(\Phi^{-{t_0}}(z)).
\end{equation}
Hence, $r_z^{-t}(I) \subset D_{\delta}(\Phi^{-t}(z))$
 for any $t$ sufficiently close to $t_0$.
It contradicts the choice of $t_0$.
\end{proof}

\begin{prop}
\label{prop:stable manifolds 2}
There exists $\epsilon_2>0$ such that
\begin{equation*}
 \bigcap_{t \geq 0}r_z^t(D_\epsilon(\Phi^{-t}(z)))
  \subset W^{cu}_\epsilon(z)
\end{equation*}
 for any $z \in \Lambda$ and $0<\epsilon<\epsilon_2$.
\end{prop}
\begin{proof}
Let $\epsilon_1>0$ be the constant 
 obtained by applying Proposition \ref{prop:stable manifolds 1}
 for $\delta=\Delta_2$.
By the uniform transversality of $(\Delta_1,E^s)$-interval to $E^s$,
 we can take a constant $\epsilon_2 \in (0,\epsilon_1)$ which satisfies
 the following property:
 For any $z \in \Lambda$ and
 $z' \in D_{\epsilon_2}(z) \- W^{cu}_{\epsilon_1}(z)$,
 there exists an $(E^u,1)$-transversal $J$ in $D_{\Delta_2}(z)$
 and $z_J \in J \cap W^{cu}_{\epsilon_1}(z)$ such that
 $\del J=\{z',z_J\}$.

Suppose that the proposition does not hold.
Then, there exists $z \in \Lambda$
 and $z' \in \bigcap_{t \geq 0}r_z^t(D_{\epsilon_2}(\Phi^{-t}(z)))
 \setminus W^{cu}_{\epsilon_2}(z)$.
Take an $(E^u,1)$-transversal $J$ in $D_{\Delta_2}(z)$
 and $z_J \in J \cap W^{cu}_{\epsilon_1}(z)$ such that
 $\del J=\{z',z_J\}$.
Since both $r_z^{-t}(z')$ and $r_z^{-t}(z_J)$ are contained in
 $D_{\Delta_2}(\Phi^{-t}(z))$ for any $t \geq 0$,
 $J$ is an $(\Delta_1/2,E^u)$-interval
 by Lemma \ref{lemma:E-u interval}.
Hence, we have
\begin{equation*}
r_z^{-t}(J) \subset D_{(\Delta_1/2)+{\Delta_2}}(\Phi^{-t}(z))
 \subset D_{\Delta_1}(\Phi^{-t}(z)). 
\end{equation*}
By Lemma \ref{lemma:Denjoy},
 the set $\Sigma=\cl{\bigcup_{t \geq 0}r_z^{-t}(J)}$ intersects with
 $\Omega_*^s \cup \Per_\irr^s(\Phi)$.
However, $\Sigma$ is contained
 in the $\Delta_1$-neighborhood of $\Lambda$.
It contradicts the choice of $\Delta_1$.
\end{proof}
We define a family $\{V^u(z)\}_{z \in \Lambda}$ of
 subsets of $M$ by
\begin{equation*}
V^u(z)=\bigcup_{t>0}\bigcup_{z' \in \cO(z)} \Phi^t(W^{cu}_\epsilon(z')).
\end{equation*}
It is a continuous family
 of $C^2$ open immersed surfaces tangent to $E^u$
 by Corollary \ref{cor:W-cu tangent}.
By Proposition \ref{prop:stable manifolds 1},
 $V^u(z)$ does not depend on the choice of
 sufficiently small $\epsilon>0$.
It is easy to see that
\begin{itemize}
\item $V^u(z_0)$ is diffeomorphic to $S^1 \times \RR$
 for any periodic point $z_0 \in \Lambda$, and 
\item  $V^u(z_1) \cap V^u(z_2) \neq \emptyset$
 for $z_1,z_2 \in \Lambda$ implies $V^u(z_1)=V^u(z_2)$.
\end{itemize}

Similar to $\{W^{cu}_\delta(z)\}_{z \in \Lambda}$,
 we can take a family $\{W^{cs}_\delta(z)\}_{z \in \Lambda}$ of
 $(E^u,1)$-transversals.
We define a family $\{V^s(z)\}_{z \in \Lambda}$ by
\begin{displaymath}
V^s(z)=\bigcup_{t>0}\bigcup_{z' \in \cO(z)} \Phi^{-t}(W^{cs}_\epsilon(z'))
\end{displaymath}
 for any small $\epsilon>0$.
It has analogous properties to $\{V^u(z)\}_{z \in \Lambda}$.

For $\Phi$-invariant compact subsets $\Lambda_1$ and $\Lambda_2$
 of $\Lambda$, we write $\Lambda_1 \preceq \Lambda_2$
 if $W^s(\Lambda_1) \cap W^u(\Lambda_2) \neq \emptyset$.
\begin{prop}
\label{prop:spectral decomposition}
Suppose that $\Lambda$ is locally maximal,
 {\it i.e.}, there exists a neighborhood $U$ of $\Lambda$
 such that $\Lambda=\bigcap_{t \in \RR} \Phi^t(U)$.
Then,
\begin{itemize}
 \item $W^s(\Lambda)=\bigcup_{z \in \Lambda \cap \Omega(\Phi)} V^s(z)$,
 \item there exists a decomposition
 $\Lambda \cap \Omega(\Phi)=\bigcup_{i=1}^m \Lambda_i$
 into mutually disjoint topologically transitive compact invariant subsets,
 and
 \item $\preceq$ is a partial order on $\{\Lambda_1,\cdots \Lambda_m\}$.
\end{itemize}
\end{prop}
\begin{proof}
By Propositions \ref{prop:stable manifolds 1}
 and \ref{prop:stable manifolds 2},
 $\Lambda$ has the shadowing property (see {\it e.g.} \cite{Sh}).
We can show the required properties
 by the same argument as the case of locally maximal hyperbolic sets.
\end{proof}
The partially ordered set $(\{\Lambda_1,\cdots, \Lambda_m\}, \preceq)$
 is called {\it the spectral decomposition} of
 $\Lambda \cap \Omega(\Phi)$.

We say a point $z$ of a topological space $X$
 is {\it accessible} from a subset $A$ of $X$
 if there exists a continuous map $l:[0,1] \ra X$ such that
 $l(1)=z$ and $l(t) \in A$ for any $t \in [0,1)$.
\begin{lemm}
\label{lemma:boundary point}
Let $\Lambda'$ be a topologically transitive
 compact invariant subset of $\Lambda$
 such that $W^s(\Lambda') \cap W^u(\Lambda')=\Lambda'$.
If $z \in \Lambda'$ is accessible from $V^s(z) \-\Lambda'$,
 then $V^u(z)$ contains a periodic point $z_* \in \Lambda'$
 which is accessible from $V^s(z_*) \- \Lambda'$.
Similarly, if $z \in \Lambda'$ is accessible from $V^u(z) \-\Lambda'$,
 then $V^s(z)$ contains a periodic orbit $z_*$
 which is accessible from $V^u(z_*) \- \Lambda'$.
\end{lemm}
\begin{proof}
The same argument as the proof of Proposition 1 of \cite{NP}
 shows that $V^u(z)$ contains a periodic point $z_*$.
Since the $\alpha$-limit set of $z$ coincides with the orbit of $z_*$,
 the invariance of $\Lambda'$ implies that $z_* \in \Lambda'$.
By the accessibility of $z$ from $V^s(z) \- \Lambda'$,
 there exists a curve $I \subset V^s(z)$ transverse to $E^s$
  such that $I \cap \Lambda'=\{z\}$.
Suppose that $z_*$ is not accessible from $V^s(z_*) \- \Lambda'$.
By the continuity of $V^u(z')$ with respect to $z' \in \Lambda$
 there exists $z_1 \in \Lambda' \cap V^u(z_*)$
 such that $V^u(z_1) \cap (I \- \{z\}) \neq \emptyset$.
Since $W^s(\Lambda') \cap W^u(\Lambda')=\Lambda'$
 and $V^\sigma(z') \subset W^\sigma(\Lambda')$
 for any $z' \in \Lambda'$ and $\sigma=s,u$,
 we have $I \cap V^u(z_1) \subset \Lambda' $.
However, it contradicts the choice of $I$.
Therefore, $z_*$ is accessible from $V^s(z_*) \- \Lambda'$.

We obtain the latter from the former by reversing the time.
\end{proof}

%
%

\subsection{Dichotomy under regularity of periodic orbits}
\label{sec:reduced dichotomy}
The aim of this subsection is to show
 Proposition \ref{prop:dichotomy} under some additional assumptions.
\begin{prop}
\label{prop:reduced dichotomy}
Let $\Phi$ be a $C^2$ $\PA$ flow
 and $TM=E^s + E^u$ be its $\PA$ splitting.
Suppose that $\Per_\irr(\Phi) \subset \Omega_*$
 and $E^s$ generates a $C^2$ foliation.
Then, either
\begin{enumerate}
\item $\Omega_*=\Per_\irr(\Phi)=\emptyset$
 and $\Phi$ is topologically transitive, or
\item $M= W^u(\Omega^s_*) \cup \Omega^u_*
 = W^s(\Omega^u_*) \cup \Omega^s_*$.
\end{enumerate}
\end{prop}
Remark that we do not assume the $C^2$-regularity of $E^u$.

Put $\Omega_h=M\-(W^u(\Omega^s_*) \cup W^s(\Omega^u_*))$.
\begin{lemm}
\label{lemma:locally maximal}
$\Omega_h$ is a locally maximal closed invariant set.
\end{lemm}
\begin{proof}
Since $\Omega_*^s$ is normally repelling,
 there exists a compact neighborhood $K^s$ of $\Omega^s_*$
 such that $\Phi^{-t}(K^s) \subset K^s$ for any $t>0$,
 $\bigcap_{t \geq 0} \Phi^{-t}(K^s)=\Omega^s_*$,
 and $\bigcup_{t \geq 0}\Phi^t(K^s)=W^u(\Omega^s_*)$.
Similarly,
 there exists a compact neighborhood $K^u$ of $\Omega^u_*$
 such that $\Phi^t(K^u) \subset K^u$ for any $t >0$,
 $\bigcap_{t \geq 0} \Phi^t(K^u)=\Omega^u_*$,
 and $\bigcup_{t \geq 0}\Phi^{-t}(K^u)=W^u(\Omega^u_*)$.
Then, $U=M \- (K^u \cup K^s)$ is a neighborhood of $\Omega_h$
 such that $\bigcap_{t \in \RR}\Phi^t(U)=\Omega_h$.
\end{proof}
It is easy to see that
 $\alpha(z) \cup \omega(z) \subset \Omega^s_* \cup \Omega^u_* \cup\Omega_h$
 for any $z \in M$.
Since $\Omega^u_*$ is normally attracting
 and $\Omega^s_*$ is normally repelling,
 we have
\begin{equation*}
M=W^u(\Omega_h) \cup W^u(\Omega_*^s) \cup \Omega_*^u
 =W^s(\Omega_h) \cup W^s(\Omega_*^u) \cup \Omega_*^s. 
\end{equation*}

We assume $\Omega_h \neq \emptyset$
 and show that $M=\Omega_h$ and $\Phi$ is topologically transitive.
It implies that $\Omega_*^u=\Omega_*^s=\emptyset$, and hence,
 $\Per_\irr(\Phi)=\emptyset$ by the assumption.

By $\cG(z)$, we  denote the leaf of a foliation $\cG$
 that contains a point $z$.
Let $\cF^s$ be the $C^2$ foliation generated by $E^s$.
\begin{lemm}
\label{lemma:leaf}
 $\cF^s(z)=V^s(z)$ for any $z \in \Omega_h$.
\end{lemm}
\begin{proof}
Since $V^s(z')$ is tangent to $E^s$ for any $z' \in \Omega_h$,
 it is a connected open subset of $\cF^s(z')$.
Since $\cF^s(z) \subset M\- \Omega^s_*
 =W^s(\Omega_h) \cup W^s(\Omega^u_*)$,
 we have a decomposition
\begin{equation*}
\cF^s(z)=(\cF^s(z) \cap W^s(\Omega^u_*)) \cup
 \bigcup_{z' \in \Omega_h \cap \cF^s(z)} V^s(z')
\end{equation*}
 of $\cF^s(z)$ into mutually disjoint open subsets.
It implies that $V^s(z)$ coincides with $\cF^s(z)$.
\end{proof}

By Proposition \ref{prop:spectral decomposition}
 and Lemma \ref{lemma:locally maximal},
 the invariant set $\Omega_h \cap \Omega(\Phi)$ admits
 a spectral decomposition $(\{\Lambda_1,\cdots,\Lambda_m\},\preceq)$.
Take a maximal element $\Lambda_+$ with respect to $\preceq$.
By the same argument as the hyperbolic case,
 we have $W^s(\Lambda_+) \cap W^u(\Lambda_+) \subset \Omega(\Phi)$.
The maximality implies that
 $W^s(\Lambda_+) \subset \Lambda_+ \cup W^u(\Omega_*^s)$.
Since $W^u(\Lambda_+) \cap W^u(\Omega_*^u)=\emptyset$,
 we have $W^s(\Lambda_+) \cap W^u(\Lambda_+) = \Lambda_+$.

Recall that a subset $\Lambda$ of $M$ is called
 {\it a saturated set} of $\cF^s$
 if $\cF^s(z) \subset \Lambda$ for any $z \in \Lambda$.
\begin{lemm}
\label{lemma:saturated}
 $\Lambda_+$ is a closed saturated set of $\cF^s$.
\end{lemm}
\begin{proof}
We will show $W^s(\Lambda_+) \subset \Lambda_+$.
It completes the proof of the lemma
 since $\cF^s(z)=V^s(z) \subset W^s(\Lambda_+)$
 for any $z \in \Lambda_+$.

Suppose that $W^s(\Lambda_+)  \not\subset \Lambda_+$.
Then, there exists $z_* \in \Lambda_+$
 which is accessible from $V^s(z_*) \- \Lambda_+$.
Since $W^s(\Lambda_+) \cap W^u(\Lambda_+)=\Lambda_+$,
 we can apply Lemma \ref{lemma:boundary point} to $\Lambda_+$.
Hence, we may assume that $z_*$ is a periodic point.
Accessibility implies that
 a connected component $L$ of $V^s(z_*)\- \cO(z_*)$
 is a subset of $W^s(\Lambda_+) \- \Lambda_+$,
 and hence, is contained in $W^u(\Omega^s_*)$.

Take a simple closed curve $\gamma \subset L$
 which is homotopic to $\cO(z_*)$ in $\cF^s(z_*)$.
Since $z_*$ is a $u$-regular periodic point,
 the holonomy of $\cF^s$ along $\gamma$ is non-trivial.

Since $W^u(T)$ is a connected open subset of $M$
 for any  torus $T$ in $\Omega^s_*$,
 there exists a torus $T_*$ in $\Omega^s_*$ such that
 $L \subset W^u(T_*)$.
Take an embedding $\psi:\TT^2 \times [-1,1] \ra W^u(T_*)$
 such that $\psi(\TT^2 \times 0)=T_*$
 and $\psi(\TT^2 \times \{-1,1\})$ is transverse to the flow.
There exists $t>0$ such that
 $\Phi^{-t}(\gamma) \subset \psi(\TT^2 \times (-1,1))$.
Let $\cG$ be the restriction of $\cF$ to $\Im \psi$.
Since $T_*$ is the unique compact leaf of $\cG$,
 a classification theorem of $C^2$ foliation on $\TT^2 \times [0,1]$
 due to Moussu and Roussarie \cite{MR}
 implies that $T_*$ is the only leaf of $\cG$ that has non-trivial holonomy.
It contradicts that
 the holonomy of $\cF^s$ along $\Phi^{-t}(\gamma)$ is non-trivial
 but $\Phi^{-t}(\gamma)$ is not contained in $T_*$.
\end{proof}
Recall that a leaf of a codimension-one foliation
 is called {\it semi-proper} when it accumulates to itself
 from at most one side.
We also say a leaf is {\it proper} when it does not accumulate to itself
 from either sides.
\begin{lemm}
\label{lemma:semi-proper}
Let $\cG$ be a $C^2$ codimension-one foliation of a
 closed three-dimensional manifold.
Then, any semi-proper leaf of $\cG$ diffeomorphic to $S^1 \times \RR$
 is proper and it has trivial holonomy.
\end{lemm}
\begin{proof}
Let $L$ be a leaf of $\cG$ which is diffeomorphic to $S^1 \times \RR$.

By the level theory of Cantwell and Conlon \cite{CC},
 $L$ is either proper or contained in an exceptional local minimal set.
Duminy's theorem \cite{CC2}
 implies that the end set of a semi-proper leaf
 in an exceptional local minimal set must be a Cantor set.
Since the end set of $L$ consists of two points,
 the leaf $L$ is proper.
By a stability theorem of proper leaves with finite ends
 due to Cantwell and Conlon \cite[Theorem 1]{CC1},
 $L$ has trivial holonomy.
\end{proof}

Now, we prove Proposition \ref{prop:reduced dichotomy}.
Let $\Phi$ be a $\PA$ flow satisfying the assumptions of the proposition.
Then, $\cF^s$ is a $C^2$ foliation and
 all periodic points in $\Omega_h$ are $u$-regular.
Suppose that the proposition does not hold.
Let $\Lambda_+$ be the maximal set in the spectral decomposition
 of $\Omega_h \cap \Omega(\Phi)$.
By Lemma \ref{lemma:saturated}, it is a closed saturated set of $\cF^s$.
Since the restriction of $\Phi$ to $\Lambda_+$
 is topologically transitive,
 the assumption implies $\Lambda_+ \neq M$.
In particular, $\Lambda_+$ contains a semi-proper leaf $L$ of $\cF^s$.
By Lemma \ref{lemma:boundary point},
 $L$ contains a periodic point $q$ in $\Lambda_+$, and hence,
 it is diffeomorphic to $S^1 \times \RR$.
Lemma \ref{lemma:semi-proper} implies that
 the holonomy of $\cF^s$ along the orbit of $q$ is trivial.
In particular, $q$ is a $u$-irregular periodic point.
However, it contradicts that
 all periodic points in $\Omega_h$ are $u$-regular.

%
%
\subsection{Local dynamics at periodic points}
\label{sec:local dynamics}
Let $\Phi$ be a $C^2$ $\PA$ flow
 with a $\PA$ splitting $TM=E^u + E^s$.
In this subsection, we suppose that $E^u$ and $E^s$ generate
 $C^2$ foliations $\cF^u$ and $\cF^u$, respectively.
Remark that $\Omega^\rho_*$ is a union of closed leaves of $\cF^\rho$
 for $\rho=s,u$.

The main aim of this subsection is to show the following proposition,
 which completes the proof of Proposition \ref{prop:dichotomy}
 by combining with Proposition \ref{prop:reduced dichotomy}.
\begin{prop}
\label{prop:regular point}
 $\Per^u_\irr(\Phi) \subset \Omega^u_*$
 and $\Per^s_\irr(\Phi) \subset \Omega^s_*$.
\end{prop}

Fix a family $\{\phi_z:[-1,1]^2 \ra M\}_{z \in M}$
 of $C^2$ local cross sections so that
 $\phi_z(0,0)=z$,
 $\phi_z([-1,1] \times y)$ is tangent to $E^s$, and
 $\phi_z(x \times [-1,1])$ is tangent to $E^u$
 for any $(x,y) \in [-1,1]^2$.
Let $\{r_z^t\}$ be the family of local returns associated to
 $\{\phi_z\}_{z \in M}$.

Recall that $D_\delta(z)$ be the $\delta$-ball in $\Im \phi_z$
 centered at $z$.
Let $\Delta>0$ be the constant obtained in Lemma \ref{lemma:Denjoy}.
For $0<\delta<\Delta$, put
\begin{align*}
I^s_\delta(z) & = D_\delta(z) \cap \phi_z([-1,1] \times 0),\\
I^u_\delta(z) & = D_\delta(z) \cap \phi_z(0 \times [-1,1]).
\end{align*}
By replacing $\Delta$ with a smaller one,
 we may assume that $I^u_\delta(z)$ and $I^s_\delta(z)$ are intervals
 for any $z \in M$ and $0<\delta <\Delta$.
\begin{lemm}
\label{lemma:Denjoy 2} 
Suppose that sequences $(z_n \in M)_{n \geq 1}$,
 $(\delta_n>0)_{n \geq 1}$, and $(t_n>0)_{n \geq 1}$
 satisfy the following properties:
\begin{itemize}
\item $\lim_{n \ra \infty}\delta_n=0$.
\item $r_{z_n}^t(I^s_{\delta_n}(z_n))$ is well-defined
 for any $n \geq 1$ and $0 \leq t \leq t_n$.
\item $\limsup_{n \ra \infty}|r_{z_n}^{t_n}(I^s_{\delta_n}(z_n))|>0$.
\end{itemize}
Then, any accumulation point of $\{z_n\}_{n \geq 1}$
 is contained in  $\Per_\irr^s(\Phi) \cup \Omega_*^s$.
\end{lemm}
\begin{proof}
Take an accumulation point $z_*$ of $(z_n)_{n \geq 1}$.
By taking subsequences if it is necessary,
 we may assume that $z_n$ converges to $z_*$,
 $\Phi^{t_n}(z_n)$ converges to a point $z_\infty$,
 and $r_{z_n}^{t_n}(I^s_{\delta_n}(z_n))$ converges to
 an interval $I_\infty \subset I^s_\Delta(z_\infty)$.
Remark that $t_n$ goes to infinity.
In fact, $|r_{z_n}^T(I^s_{\delta_n}(z_n))|$ converges
 to zero for any given $T>0$ since $\delta_n$ goes to zero.

The interval $I_\infty$ is a $C^2$ $(\Delta, E^u)$-interval,
By Lemma \ref{lemma:Denjoy},
 there exists $z' \in \Int I_\infty$ such that
 its $\alpha$-limit set $\alpha(z')$ is a periodic orbit
 in $\Per_\irr^s(\Phi)$ or an embedded torus in $\Omega_*^s$.
In each case, $N\Phi^{-t}|_{E^u/T\Phi}$ is uniformly contracting
 on $\alpha(z')$.
Hence,
 there exists a compact neighborhood $V$ of $z'$ in $\cF^u(z')$
 such that
 $\bigcap_{t>0}\cl{\bigcup_{t'>t}\Phi^{-t'}(V)}=\alpha(z')$.
For any sufficiently large $n \geq 1$,
 the interval $r_{z_n}^{t_n}(I^s_{\delta_n}(z_n))$
 contains a point $z'_n$ of $V$.
Then, $z_*=\lim_{n \ra \infty} (r_{z_n}^{t_n})^{-1}(z'_n)$
 is contained in  $\alpha(z')$.
Hence, $z_*$ is contained in $\Per_\irr^s(\Phi)$ or $\Omega_*^s$.
\end{proof}

\begin{lemm}
\label{lemma:open}
Let $z_0$ be an $s$-regular periodic point.
Then, the following holds:
\begin{itemize}
\item  There exists $\delta>0$
 and $\tau: W^s(\cO(z_0)) \ra \RR$
 such that $I^s_\delta(\Phi^{\tau(z)}(z)) \subset W^s(\cO(z_0))$
 for any $z \in W^s(\cO(z_0))$.
\item For any $z \in W^s(\cO(z_0))$,
 $\cF^s(z) \cap W^s(\cO(z_0))$ is an open subset
 of $\cF^s(z)$ with respect to the leafwise topology.
\item If $A$ is any open annulus
 such that $\cO(z_0) \subset A \subset \cF^s(z_0) \cap W^s(\cO(z_0))$,
 then $\bigcup_{t \geq 0} \Phi^{-t}(A)$ is a connected component
 of $\cF^s(\cO(z_0)) \cap W^s(\cO(z_0))$.
\end{itemize}
\end{lemm}
\begin{proof}
Let $T$ be the period of $z_0$.
There exists a closed interval $I \subset [-1,1]$
 and $C^2$ maps $f,g:I \ra [-1,1]$ 
 such that $0 \in f(I) \subset \Int I$,
 $r_{z_0}^T \circ \phi_{z_0}(x,y)=\phi_{z_0}(f(x),g(y))$,
 and $\bigcap_{n \geq 0} f^n(I)=\{0\}$.
Put $\Lambda^u=\bigcap_{n \geq 0}g^{-n}(I)$ and
$\Lambda^u_0=\{y \in \Lambda^u \st \lim_{n \ra \infty} g^n(y)=0 \}$.
For any $(x,y) \in \Int I \times \Lambda^u$,
 we have
 $r_{z_0}^{nT} \circ \phi_{z_0}(x,y)=\phi_{z_0}(f^n(x),g^n(y))$.
By the compactness of $\Lambda^u$,
 there exists $\delta>0$ such that
\begin{equation}
\label{eqn:open 1} 
I^s_\delta(\phi_{z_0}(x,y)) \subset \phi_{z_0}(\Int I \times y)
\end{equation}
 for any $(x,y) \in f(I) \times \Lambda^u$.

If $(x,y) \in \Int I \times \Lambda^u_0$, then
 $(f^n(x),g^n(y))$ converges to $(0,0)$ as $n$ goes to infinity.
For any $z \in W^s(\cO(z_0))$,
 its positive orbit $\{\Phi^t(z) \st t  \geq 0\}$
 intersects with $\phi_{z_0}(f(I) \times \Lambda^u_0)$.
Hence, we have
\begin{equation*}
 W^s(\cO(z_0))=\bigcup_{t \geq 0}
 \Phi^{-t} \circ \phi_{z_0}(f(I) \times \Lambda^u_0)
 =\bigcup_{t \geq 0}
 \Phi^{-t} \circ \phi_{z_0}(\Int I \times \Lambda^u_0).
\end{equation*}
It completes the proof of the first assertion of the lemma
 with the inclusion (\ref{eqn:open 1}).
The second assertion is an immediate consequence of the first.

Put
\begin{align*}
V_1 & = \bigcup_{t \geq 0} \Phi^{-t} \circ \phi_{z_0}(\Int I \times 0),\\
V_2 & = \bigcup_{t \geq 0} \Phi^{-t} \circ \phi_{z_0}
 (\Int I \times (\Lambda^u_0 \- \{0\})).
\end{align*}
Since $f(I) \subset \Int I$,
 $g(0)=0$, and $g(\Lambda^u_0\-\{0\}) \subset \Lambda^u_0\-\{0\}$,
 the set $\cF^s(z_0) \cap W^s(\cO(z_0))$
 is a disjoint union of its open subsets
 $V_1$ and $V_2 \cap \cF^s(z_0)$.
Hence, $V_1$ is a connected component of $\cF^s(z_0) \cap W^s(\cO(z_0))$.
Let $A$ be the annulus in the third assertion of the lemma.
Then, $\phi_{z_0}(f^m(I)) \subset A$ for some large $m \geq 1$.
It implies
\begin{equation*}
V_1=\bigcup_{t \geq 0} \Phi^{-t} \circ \phi_{z_0}(f^m(I) \times 0)
 \subset \bigcup_{t \geq 0}\Phi^{-t}(A) \subset V_1.
\end{equation*}
\end{proof}

\begin{lemm}
\label{lemma:annulus}
Let $z_0$ and $z_1$ be periodic points of $\Phi$.
\begin{enumerate}
\item If $z_0$ is $s$-regular and
 $z_1$ is accessible from
 a connected component $V$ of $W^s(\cO(z_0)) \cap \cF^s(z_1)$.
 then $\cF^s(z_0)=\cF^s(z_1)$,
 the orbits of $z_0$ and $z_1$ are homotopic in $\cF^s(z_0)$
 as unoriented curves,
 and $V$ contains $\cO(z_0)$.
\item If $z_0$ attracting and
 $z_1$ is accessible from $W^s(\cO(z_0)) \cap \cF^u(z_1)$,
 then $\cF^u(z_0)=\cF^u(z_1)$ and
  the orbits of $z_0$ and $z_1$ are homotopic in $\cF^u(z_0)$
 as unoriented curves.
\end{enumerate}
\end{lemm}
\begin{proof}
First, we show the former assertion of the lemma.
Suppose that $z_0$ is $s$-regular.
Let $T$ be the period of $z_0$.
There exists a closed interval $I$
 and $C^2$ maps $f,g: I \ra [-1,1]$ such that
 $0 \in \Int I \subset f(I)$,
 $\bigcap_{n \geq 0} f^n(I)=\{0\}$,
 $\phi_{z_0}(I \times I) \cap \cO(z_0)=\{z_0\}$, and 
 $r_{z_0}^T \circ \phi_{z_0}(x,y)=\phi_{z_0}(f(x),g(y))$
 for any $(x,y) \in I \times I$.
Put $U=\bigcup_{t=0}^T r_{z_0}^t(I \times I)$
 and let $\cG(y)$ be the connected component
 of $\cF^s(\phi_{z_0}(0,y)) \cap U$ which contains $\phi_{z_0}(0,y)$.
Then, we can see the following properties of $\cG$ and $U$:
\begin{itemize}
\item $\cG(y)$ is not contractible if and only if $g(y)=y$.
\item If $z=\phi_{z_0}(x,y)$ is a point of $W^s(\cO(z_0))$
 and $\Phi^t(z) \in U$ for any $t \geq 0$,
 then $g^n(y)$ converges to $0$ as $n$ tends to infinity.
\end{itemize}

Suppose that a periodic point $z_1$ is accessible from
 a connected component $V$ of $\cF^s(z_1) \cap W^s(\cO(z_0))$.
Then, there exists a simple closed curve $\gamma$ in $V$
 which is homotopic to $\cO(z_1)$.
By the Poincar\'e-Bendixon theorem,
 $\cO(z_1)$ and $\gamma$ are not null-homotopic in $\cF^s(z_1)$.
We can take $t_1>0$ such that $\Phi^t(\gamma) \subset U$
 for any $t \geq t_1$.
It implies that $\Phi^{t_1}(\gamma) \subset \cG(y)$
 for some $y \in \bigcap_{n \geq 0} g^{-n}(I)$
 with $\lim_{n \ra \infty}g^n(y)=0$.
Since the closed curve $\Phi^{t_1}(\gamma)$ is not null-homotopic in
 $\cG(y) \subset \cF^s(z_1)$,
 we have $y=0$.
Hence, $\Phi^{t_1}(\gamma)$ is homotopic
 to $\cO(z_0)$ in $\cF^s(z_0)$ as an unoriented curve
 and the set $V$ intersects the connected component of
 $\cF^s(z_0) \cap W^s(\cO(z_0))$
 which contains $\cO(z_0)$.
Therefore, we have $\cF^s(z_1)=\cF^s(z_0)$,
 $\cO(z_0)$ and $\cO(z_1)$ are homotopic in $\cF^s(z_0)$
 as unoriented curves, and $\cO(z_0) \subset V$.

Next, we show the latter assertion of the lemma.
Suppose that $z_0$ is attracting.
Let $T$ be the period of $z_0$.
There exists a closed interval $I$
 and $C^2$ maps $f_1,g_1: I' \ra [-1,1]$ such that
 $0 \in \Int I' \subset g_1(I')$,
 $\bigcap_{n \geq 0} g_1^n(I')=\{0\}$,
 $0$ is the unique fixed point of $f_1|_{I'}$,
 $\phi_{z_0}(I \times I') \cap \cO(z_0)=\{z_0\}$, and 
 $r_{z_0}^T \circ \phi_{z_0}(x,y)=\phi_{z_0}(f_1(x),g_1(y))$
 for any $(x,y) \in I' \times I'$.
Put $U'=\bigcup_{t=0}^T r_{z_0}^t(I \times I)$
 and let $\cG'(x)$ be the connected component
 of $\cF^u(\phi_{z_0}(x,0)) \cap U'$ which contains $\phi_{z_0}(x,0)$.
Then, we can see the following properties of $\cG'$ and $U'$:
\begin{itemize}
\item $\cG'(x)$ is not contractible if and only if $x=0$.
\item If $z=\phi_{z_0}(x,y)$ is a point of $W^s(\cO(z_0))$
 and $\Phi^t(z) \in U'$ for any $t \geq 0$,
 then $f_1^n(x)$ converges to $0$ as $n$ tends to infinity.
\end{itemize}
Now, the same argument as above, where we replace $g$, $\cG$, and $\cF^s$
 with $f_1$, $\cG'$ and $\cF^u$, respectively,
 shows the latter assertion of the lemma.
\end{proof}

\begin{lemm}
\label{lemma:stable set}
The following holds for any $s$-regular periodic point $z_0$:
\begin{enumerate}
 \item $\cF^s(z) \subset W^s(\cO(z_0))$
 for any $z \in W^s(\cO(z_0)) \- \cF^s(z_0)$.
 \item $\cF^s(z_0) \cap W^s(\cO(z_0))$
 is homeomorphic to $S^1 \times \RR$.
 \item If $\cF^s(z_0) \not\subset W^s(\cO(z_0))$,
 then $\cF^s(z_0)$ contains an $s$-irregular periodic point
 such that the orbits of $z_0$ and $z_1$ are homotopic 
 as unoriented closed curves in $\cF^s(z_0)$.
\end{enumerate}
\end{lemm}
\begin{proof}
Since $z_0$ is $s$-regular,
 there exists an embedded closed annulus $A_0 \subset \cF^s(z_0)$
 such that $\Phi^t(A_0) \subset \Int A_0$ for any $t>0$
 and $\bigcap_{t >0}\Phi^t(A_0)=\cO(z_0)$.
Put $V_0=\bigcup_{t \geq 0} \Phi^{-t}(A_0)$.
It is diffeomorphic to $S^1 \times \RR$
 and is a connected component of $W^s(\cO(z_0)) \cap \cF^s(z_0)$
 by the third item of Lemma \ref{lemma:open}.

Fix a leaf $L$ of $\cF^s$
 and a connected component $V$ of $L \cap W^s(\cO(z_0))$.
We suppose that $V \neq L$ and show $V=V_0$.
Take $z_1 \in L \- V$ which is accessible from $V$.
There exist sequences $(z'_n \in V)_{n \geq 1}$
 and $(\delta_n>0)_{n \geq 1}$ such that
 $z_1 \in I^s_{\delta_n}(z'_n)$ for any $n \geq 1$
 and $\delta_n$ converges to zero.
By Lemma \ref{lemma:open},
 we can choose $\delta>0$ and $(T_n>0)_{n \geq 1}$
 such that $I^s_{\delta}(\Phi^{T_n}(z'_n)) \subset W^s(\cO(z_0))$
 for any $n \geq 1$.
Since $\cO(z_1) \cap W^s(\cO(z_0))= \emptyset$,
 there exists a sequence $(t_n \in (0,T_n))_{n \geq 1}$ such that
 $r_{z_n'}^t(I^s_{\delta_n}(z'_n))$ is well-defined for any $t \in [0,t_n]$
 and $|r_{z_n'}^{t_n}(I^s_{\delta_n}(z'_n))|=\delta$.
By Lemma \ref{lemma:Denjoy 2},
 $z_1$ is a point of $\Per_\irr^s(\Phi) \cup \Omega_*^s$.
Since $\cF^s(z_1)$ contains a periodic point $z_1$
 and some non-periodic points in $V$,
 we have $\cF^s(z_1) \not\subset \Omega_*^s$.
Therefore, $z_1$  is an $s$-irregular periodic point.
Since $z_1$ is accessible from $V \subset W^s(\cO(z_0))$,
 we can apply the former part of Lemma \ref{lemma:annulus}.
It implies that $\cF^s(z_1)=\cF^s(z_0)$,
 $V=V_0$, and the orbits of $z_0$ and $z_1$ are homotopic
 as unoriented closed curves in $\cF^s(z_0)$.
\end{proof}

\begin{lemm}
\label{lemma:trivial holonomy}
$\Per^s_\irr(\Phi) \cap \cF^s(z)$ is
 a closed subset of $\cF^s(z)$ for any $z \in M$.
Similarly, $\Per^u_\irr(\Phi) \cap \cF^u(z)$
 is a closed subset of $\cF^u(z)$ for any $z \in M$.
\end{lemm}
\begin{proof}
We show the former assertion.
The proof of the latter is similar.

Suppose that a sequence $(z_n)_{n \geq 1}$
 in $\cF^s(z) \cap \Per_\irr^s(\Phi)$ converges to
 $z_*$ with respect to the leafwise topology of $\cF^s(z)$.
By replacing $z_n$ in its orbit,
 we may assume that there exists a sequence $(\delta_n>0)_{n \geq 1}$
 such that $\lim_{n \ra \infty}\delta_n=0$
 and $I^s_{\delta_n}(z_n)$ contains $z_*$ for any $n \geq 1$.
Let $T_n$ be the period of $z_n$ 
 and $\Delta>0$ be a constant
 such that $D_\Delta(z) \subset \Im \phi_z$ for any $z \in M$.

First, we suppose that
 there exists $n_0 \geq 1$ such that
 $|r_{z_{n_0}}^t(I^s_{\delta_{n_0}}(z_{n_0}))| \leq \Delta$
 for any $t \in [0,T_{n_0}]$.
Then, $r_{z_{n_0}}^{T_{n_0}}(I^s_{\delta_{n_0}}(z_{n_0}))$
 is well-defined.
Since $z_*$ is a limit point of the sequence $(z_n)_{n \geq 1}$
 of periodic points with respect to the leafwise topology,
 we can choose a sequence $(z_n')_{n \geq 1}$
 of periodic points of $\Phi$ such that
 $\lim_{n \ra \infty}z_n'=z_*$
 and $z_n' \in I^s_{\delta_{n_0}}(z_{n_0})$ for any $n$.
Since $r_{z_{n_0}}^{T_{n_0}}(z_{n_0})=z_{n_0}$
 and $r_{z_{n_0}}^{T_{n_0}}$ is an orientation preserving
 homeomorphism on an interval,
 we have $r_{z_{n_0}}^{T_{n_0}}(z_n')=z_n'$ for any $n$.
It implies that $r_{z_{n_0}}^{T_{n_0}}(z_*)=z_*$.
Therefore, $z_*$ is an $s$-irregular periodic point of $\Phi$.

Now, it is sufficient to consider the case that
 there exists a sequence $(t_n \in [0,T_n])_{n \geq 1}$ such that
 $|r_{z_n}^t(I^s_{\delta_n}(z_n))|\leq \Delta$
 for any $t \in [0,t_n]$ and
 $|r_{z_n}^{t_n}(I^s_{\delta_n}(z_n))|=\Delta$.
By Lemma \ref{lemma:Denjoy 2},
 $z_*=\lim_{n \ra \infty}z_n$ is a point of
 $\Per_{\irr}^s(\Phi) \cup \Omega_*^s$.
If $z_* \in \Omega_*^s$, then
 $\cF^s(z_*)$ is an embedded torus contained in $\Per_\irr^s(\Phi)$
 since $\{z_n\} \subset \Per_\irr^s(\Phi) \cap \cF^s(z_*)$.
\end{proof}

\begin{lemm}
\label{lemma:trivial holonomy 2} 
For any $z_0 \in \Per_\irr^u(\Phi)$,
 the leaf $\cF^s(z_0)$ is proper, has trivial holonomy,
 and is contained in $W^s(\cO(z_0))$.
\end{lemm}
\begin{proof}
The periodic point $z_0$ is $s$-regular.
Hence, there exists a closed annulus $A^s$
 in $\cF^s(z_0) \cap W^s(\cO(z_0))$
 such that $\Phi^t(A^s) \subset \Int A^s$ for any $t>0$
 and $\bigcap_{t \geq 0}\Phi^t(A^s)=\cO(z_0)$.
Put $V=\bigcup_{t \geq 0}\Phi^{-t}(A^s)$.
Lemmas \ref{lemma:open} and \ref{lemma:stable set}
 imply that $V=\cF^s(z_0) \cap W^s(\cO(z_0))$
 and $\cF^s(z) \subset W^s(\cO(z_0))$
 for any $z \in W^s(\cO(z_0)) \- V$.
Since $z_0$ is $u$-irregular,
 there exists a closed annulus $A^u_0$
 in $\cF^u(z_0)$ such that $\cO(z_0) \subset \del A^u_0$,
 $A^u_0 \cap A^s=\cO(z_0)$,
 and $\Phi^t(A^u_0) \subset A^u_0$ for any $t \geq 0$.
It is easy to see that $V \cap A^u_0=\cO(z_0)$.

First, we show $V=\cF^s(z_0)$.
Suppose that it does not hold.
By Lemma \ref{lemma:stable set},
 $\cF^s(z_0)$ contains an $s$-irregular periodic point $z_1$
 such that the orbits of $z_0$ and $z_1$ are homotopic
 as unoriented closed curves in $\cF^s(z_0)$.
For $i=0,1$, let $T_i$ be the period of $z_i$.
Put $\lambda^u_i=\|D\Phi^{T_i}|_{E^u(z_i)}\|$
 and $\lambda^s_i=\|D\Phi^{T_i}|_{E^s(z_i)}\|$.
Since $z_0$ is $u$-irregular, $z_1$ is $s$-irregular,
 and $TM=E^s+ E^u$ is a dominated splitting, we have
\begin{equation}
\label{eqn:trivial holonomy 2-1} 
\lambda_0^s<\lambda_0^u \leq 1 \leq \lambda_1^s < \lambda_1^u.
\end{equation}
Since $\lambda^u_i$ is the absolute value of
 the linear holonomy of $\cF^s$ along the orbit of $z_i$ for each $i$
 and the orbits of $z_0$ and $z_1$ are homotopic in $\cF^s(z_0)$
 as unoriented curves in $\cF^s(z_0)$,
 we have $\lambda_0^u=(\lambda_1^u)^{\pm 1}$.
Hence, the inequality (\ref{eqn:trivial holonomy 2-1}) implies
\begin{equation}
\label{eqn:trivial holonomy 2-2} 
\lambda_0^u=(\lambda_1^u)^{-1} <1. 
\end{equation}
In particular, $\cO(z_0)$ is an attracting periodic orbit.
Since $W^s(\cO(z_0))$ is an open subset of $M$,
 there exists a closed annulus $A^u_1 \subset \cF^u(z_1)$
 such that $\cO(z_1) \subset \del A^u_1$
 and $\cF^s(z)$ intersects with $(A^u_0\- \cO(z_0)) \cap W^s(\cO(z_0))$
 for any $z \in A^u_1 \- \cO(z_1)$.
Now, we recall that $V \cap A^u_0=\cO(z_0)$
 and $\cF^s(z) \subset W^s(\cO(z_0))$ for any $z \in W^s(\cO(z_0)) \- V$.
They imply that $A^u_1 \- \cO(z_1)$ is contained in $W^s(\cO(z_0))$.
In particular, $z_1$ is accessible from $\cF^u(z_1) \cap W^s(\cO(z_0))$.
By Lemma \ref{lemma:annulus}, we have $\cF^u(z_1)=\cF^u(z_0)$
 and the orbits of $z_0$ and $z_1$ are homotopic in $\cF^u(z_0)$
 as unoriented closed curves.
Since $\lambda^s_i$ is the absolute value of
 the linear holonomy of $\cF^u$ along the orbit of $z_i$,
 we have $\lambda^s_0=(\lambda^s_1)^{\pm 1}$.
Hence, the inequality (\ref{eqn:trivial holonomy 2-1}) implies
\begin{equation*}
\lambda_0^s=(\lambda_1^s)^{-1}.
\end{equation*}
However, it contradicts with the inequalities
 (\ref{eqn:trivial holonomy 2-1}) and (\ref{eqn:trivial holonomy 2-2}).
Therefore, we have $V=\cF^s(z_0)$.

Since $V \cap A^u_0=\cO(z_0)$, the leaf $\cF^s(z_0)=V$ is semi-proper.
By Lemma \ref{lemma:semi-proper},
 $\cF^s(z_0)$ is a proper leaf with trivial holonomy.
\end{proof}

Now, we prove Proposition \ref{prop:regular point}.
We claim that $\Per^u_\irr(\Phi) \subset \Omega^u_*$.
Once it is done, we can show that
 $\Per^s_\irr(\Phi) \subset \Omega^s_*$
 by applying the claim to the time-reverse of $\Phi$.

Take a leaf $L^u$ of $\cF^u$ which intersects with $\Per_\irr^u(\Phi)$.
By Lemmas \ref{lemma:trivial holonomy} and \ref{lemma:trivial holonomy 2},
 $\Per_\irr^u(\Phi) \cap L^u$ is a closed and open subset of $L^u$.
Hence, $L^u$ is a subset of $\Per_\irr^u(\Phi)$.
It is sufficient to show that $L^u$ is a closed leaf.
If it is not, then there exists a transversal $J$ of $\cF^u$
 such that $J$ is contained in a leaf $L^s$ of $\cF^s$
 and $J \cap L^u$ is not a relatively compact subset of $L^u$.
Take a point $z_* \in L^u \cap J$.
By Lemma \ref{lemma:trivial holonomy 2},
 $J \subset L^s$ is contained in $W^s(\cO(z_*))$.
Since $L^u \subset \Per(\Phi)$,
 it implies that $L^u \cap J$ is contained in $\cO(z_*)$.
However, it contradicts that $L^u \cap J$ is not a relatively compact
 subset of $L^u$.

\section{Topologically transitive regular $\PA$ flows}
\label{sec:transitive}

%
%
\label{sec:section3}
The aim of this section is the following proposition,
 which completes the proof of the main theorem
 by combining with Propositions \ref{prop:dichotomy}
 and \ref{prop:T2*I}.
\begin{prop}
\label{prop:Anosov} 
If a $C^2$ topologically transitive $\PA$ flow
 on a closed three-dimensional manifold
 admits a $C^2$ $\PA$ splitting, then it is an Anosov flow.
\end{prop}

When the $\PA$ flow $\Phi$ admits a global cross section,
 the proof of Proposition \ref{prop:Anosov}
 is reduced to an observation that the distortion
 of a holonomy map of a one-dimensional foliation on a surface
 can be estimated by the area of rectangle sweeped out by the holonomy.
We refer the reader to \cite{As2} for the detail.
If a $\PA$ flow admits invariant one-dimensional subbundles
 of $E^s$ and $E^u$ which are transverse to the flow,
 then we can apply the proof in \cite{As2} with a small modification.
However, a $\PA$ flow admits no invariant one-dimensional
 subbundles transverse to the flow in general.
This is the main technical difficulty in the proof.

The structure of the section is as follows:
In Subsections \ref{sec:one-dim} and \ref{sec:Markov},
 we show the (possibly non-uniform) contraction along
 the direction transverse to $E^u$
 by the standard argument using a Markov partition
 and a theorem due to Ma\~n\'e.
In Subsection \ref{sec:non-expansion},
 in order to overcome the above technical difficulty,
 we show that the diameter of $\Phi^t(D^s)$ 
 is uniformly bounded for any small disk $D^s$ tangent to $E^s$
 and $t \geq 0$
 after we replace the original flow $\Phi$ by a suitable time-change.
This condition makes us enable to apply the method in \cite{As2}.
It is done in Subsection \ref{sec:hyp}.

\subsection{One-dimensional topological Markov maps}
\label{sec:one-dim}

In this section, we prove some results about piecewise $C^2$
 Markov maps on a finite union of compact intervals,
 which we use later.

Let $I_*$ be a finite union of compact intervals in $\RR$
 and $\Lambda_*$ be a finite subset of $\Int I_*$.
We say a map $F:I_* \ra I_*$ is {\it a $C^2$ pre-Markov map}
 with {\it the set of discontinuity $\Lambda_*$} if 
 $F(\Lambda_*) \subset \del I_*$,
 $\cl{F(I_* \- \Lambda_*)}=I_*$,
 and for each connected component $J$ of $I_* \- \Lambda_*$,
 there exists a connected component $I$ of $I_*$ such that
 $F|_J$ extends to a $C^2$ diffeomorphism from $\cl{J}$ to $I$.
Put $I_*^n=I_* \-\bigcup_{n'=0}^nF^{-n'}(\Lambda_*)$
 and $I_*^\infty=\bigcap_{n \geq 0} I_*^n$.
For $x \in I_*^\infty$ and $n \geq 0$,
 let $I_*^n(x)$ be the connected component of $I_*^n$ that contains $x$.
Then, the restriction of $F^m$ to $\cl{I_*^n(x)}$
 extends to a $C^2$ diffeomorphism
 onto $\cl{I_*^{n-m}(F^m(x))}$ for any $x \in I_*^n$
 and $n \geq m \geq 0$.

For $n \geq 1$, let $\Per_n(F)$ be the set of all periodic points
 of period $n$.
Since $F(\Lambda_* \cup \del I_*) \subset \del I_*$
 and $\del I_* \cap \Lambda_*=\emptyset$,
 the set $\Per_n(F)$ is a subset of $I_*^\infty$
 for any $n \geq 1$.

We say a pre-Markov map $F$
 is {\it a $C^2$ topologically Markov map} if 
 $\bigcap_{n \geq 1}\cl{I_*^n(x)}$ consists of exactly one point
 for any $x \in I_*^\infty$.
\begin{lemm}
\label{lemma:Markov 1} 
Let $F$ be a $C^2$ topologically Markov map on $I_*$.
Then, $\Per_n(F)$ is finite for any $n \geq 1$
 and any periodic point $x \in \Per_n(F)$ satisfies
 $|DF^n(x)|\geq 1$.
\end{lemm}
\begin{proof}
It is an immediate consequence of the equations
 $\bigcap_{k \geq 0}\cl{I_*^{kn}(x)}=\{x\}$ and
\begin{equation*}
 F^n\left(\cl{I_*^{(k+1)n}(x)}\right)
  =\cl{I_*^{kn}(x)}  \supset \cl{I_*^{(k+1)n}(x)}
\end{equation*}
 for any $x \in \Per_n(F)$.
\end{proof}

Recall the definition
 and some properties of the distortion of a one-dimensional map.
For a $C^2$ map $h:I \ra I'$ between intervals $I$ and $I'$,
 we define {\it the distortion} $\dist(h,I)$ by
\begin{displaymath}
 \dist(h,I) =\sup_{x,y \in I}\left(\log|Dh(x)|-\log|Dh(y)|\right).
\end{displaymath}
It is easy to verify 
\begin{equation}
\label{eqn:distortion 2}
\dist(h_n \circ \cdots \circ h_0,I)
  \leq \sum_{m=0}^n \dist(h_m,h_{m-1} \circ \cdots \circ h_0(I)).
\end{equation}
By the Mean Value Theorem for $h$ and $\log|Dh|$, we also have
\begin{align}
\label{eqn:distortion 3}
\dist(h,I)
 &  \geq \sup_{x \in I}\left|\log |Dh(x)|
 -\log\frac{|h(I)|}{|I|}\right|,\\
\label{eqn:distortion 1}
 \dist(h,I)
 & \leq |I| \cdot \sup_{x \in I}|D(\log|Dh|)(x)|.
\end{align}

Let $\Per_*(F)$ be the set of non-hyperbolic periodic points.
The following proposition is a variant of Man\'e's theorem(\cite{Ma}).
\begin{prop}
\label{prop:Mane}
Let $F$ be a $C^2$ topologically Markov map on $I_*$
 and $\Lambda_*$ be the set of discontinuity of $F$.
Then, 
\begin{equation}
\label{eqn:Mane 1}
\lim_{n \ra \infty}
  \left( \inf\{|DF^n(x)| \st x \in \Per_n(F)\}\right)=\infty
\end{equation}
 and $\Per_*(F)$ consists of only finitely many points.
Moreover, for any given neighborhood $U$ of $\Per_*(F)$,
\begin{equation}
\label{eqn:Mane 2}
\lim_{n \ra \infty}\left(\inf\{|DF^n(x)|
 \st x \in I_*^n\-F^{-n}(U) \}\right)
 =+\infty.  
\end{equation}
\end{prop}
\begin{proof}
The equation (\ref{eqn:Mane 1})
 is a consequence of Theorem 5.1 of \cite{Ho},
 which is a version of Ma\~n\'e's theorem for piecewise $C^2$ maps. 
By Lemma \ref{lemma:Markov 1},
 it implies the finiteness of non-hyperbolic periodic points.

Since $I_*^n$ consists of finitely many connected components
 for any $n$ and $\bigcap_{n \geq 0}\cl{I_*^n(x)}=\{x\}$
 for any $x \in I_*^\infty$,
 there exist sequences $(K_n)_{n=0}^\infty$
 and $(K'_n)_{n=0}^\infty$ such that
 $\lim_{n \ra \infty}K_n=\lim_{n \ra \infty}K'_n=0$ and
 $K_n \leq |I_*^n(x)| \leq K'_n$ for any $n \geq 0$
 and $x \in I_*^n$.

Fix a neighborhood $U$ of $\Per_*(F)$.
By the finiteness of $\Per_*(F)$,
 there exists $N \geq 1$ such that
 $I_*^N(x_*) \subset U$ for any $x_* \in \Per_*(F)$.
Remark that $I_*^N(x) \cap \Per_*(F)=\emptyset$ if $x \not\in U$.
Put
\begin{align*}
 c_1
 & =|I_*| \cdot \sup\{|D(\log|DF|)|(x) \st x \in I_* \- \Lambda_*\},\\
 c_2
 & =\inf\{|DF(x)| \st x \in I_* \- \Lambda_*\}.
\end{align*}
We say that an interval $J \subset I_*$
 is $(\lambda,m)$-compatible for $\lambda>1$ and $m \geq 1$
 if $J \subset I_*^m$ and
 $\sum_{i=0}^m|F^i(J)|<\frac{\lambda}{\lambda-1} |I_*|$.
The properties of distortion imply
\begin{equation*}
\inf_{x \in J}|DF^m(x)|
 \geq  e^{-\frac{c_1\lambda}{\lambda-1}} \cdot \frac{|F^m(J)|}{|J|}
\end{equation*}
 for any $(\lambda,m)$-compatible interval $J$.

By the argument in Section III.5 of \cite{MS},
 it can be shown that there exists $\lambda>1$ and $n_0 \geq 1$
 such that $I_*^{n+N}(x)$ is a $(\lambda,n-n_0)$-compatible interval
 if $x \in I_*^{n+N}$ satisfies
 $I_*^N(F^n(x)) \cap \Per_*(F)=\emptyset$.
It implies that
\begin{align*}
 |DF^n(x)|
 & \geq |DF^{n_0}(F^{n-n_0}(x))| \cdot
   e^{-\frac{c_1 \lambda}{\lambda-1}}
  \frac{|I_*^{N+n_0}(F^{n-n_0}(x))|}{|I_*^{n+N}(x)|}\\
 & \geq
   \frac{c_2^{n_0}e^{-\frac{c_1 \lambda}{\lambda-1}} K_{N+n_0}}{ K'_n}
\end{align*}
 for any $n \geq 1$
 and $x \in I_*$ with $I_*^N(F^n(x)) \cap \Per_*(F)=\emptyset$.
Since $\lim_{n \ra \infty}K'_n=0$,
 it completes the proof. 
\end{proof}

%
%
\subsection{Markov partitions}
\label{sec:Markov}
Fix a closed three-dimensional Riemannian manifold $M$.
Let $\{(\phi_k,\tau_k)\}_{k=1}^m$ be a family
 of pairs of a continuous embedding of $[0,1]^2$ into $M$
 and a continuous positive-valued function on $[0,1]^2$
 such that the map $(w,t) \mapsto \Phi^t \circ \phi_k(w)$ is 
 an embedding of $\{(w,t) \st w \in [0,1]^2, t \in [0,\tau_k(w)]\}$
 into $M$ for each $k$.
We define a family $\{\phi'_k\}_{k=1}^m$ of embeddings
 of $[0,1]^2$ into $M$
 by $\phi'_k(w)=\Phi^{\tau_k(w)} \circ \phi_k(w)$.
Put $R_k=\Im \phi_k$, $R'_k=\Im \phi'_k$, and
 $P_k=\{\Phi^t \circ \phi_k(w) \st w \in [0,1]^2, t \in [0,\tau(w)]\}$.

We say the family $\{(\phi_k,\tau_k)\}$
 determines {\it a Markov partition} $\{P_k\}_{k=1}^m$
 associated to a flow $\Phi$
 if it satisfies the following properties:
\begin{enumerate}
\item $\phi_k([0,1] \times y)$ and $\phi_k(x \times [0,1])$
 are intervals tangent to $E^s$ and $E^u$ respectively,
 for any $(x,y) \in [0,1]^2$ and $k=1,\ldots,m$.
\item $M=\bigcup_{k=1}^m P_k$.
\item For each pair $(k,l)$,
 $\del (P_k \cap P_l)=\del P_k \cap \del P_l$
 and there exist subintervals $I_{k,l}$ and $J_{k,l}$ of $[0,1]$,
 which may be empty sets, such that
\begin{displaymath}
 R'_k \cap R_l =\phi'_k([0,1] \times I_{k,l})
 =\phi_l(J_{k,l} \times [0,1]).
\end{displaymath}
\end{enumerate}
Remark that $\{\Int I_{k,l'}\}_{l'=1}^m$ and $\{\Int J_{k',l}\}_{k'=1}^m$
 are partition of $[0,1]$ up to finite set for any fixed $(k,l)$.
Put $I_*= [0,1] \times \{1,\ldots,m\}$.
The family $\{(\phi_k,\tau_k)\}_{k=1}^m$
 induces a piecewise continuous map $F:I_* \ra I_*$ by $F(y,k)=(y',l)$
 if $\phi_l(J_{k,l} \times y')=\phi'_k([0,1] \times y)$.
We call the map $F$ {\it the reduced return map}.

We say a $\PA$ flow $\Phi$ is {\it $E^s$-fine} if
\begin{enumerate}
\item both $\Omega_*$ and $\Per_\irr(\Phi)$ are empty, and
\item $\Phi$ admits a $\PA$ splitting $TM=E^u + E^s$
 such that $E^s$ is a $C^2$ subbundle.
\end{enumerate}
Remark that $\Phi$ is topologically transitive
 by Proposition \ref{prop:reduced dichotomy}.

\begin{lemm}
\label{lemma:Markov}
Let $\Phi$ be a $C^2$ $E^s$-fine $\PA$ flow.
For any given $\epsilon>0$,
 $\Phi$ admits a Markov partition $\{P_k\}_{k=1}^m$
 such that the reduced return map is a $C^2$ topological Markov map
 and the diameter of each $P_l$ is smaller than $\epsilon$.
\end{lemm}
\begin{proof}
As we see in Section \ref{sec:dichotomy},
 the flow $\Phi$ has the shadowing property on $M$.
Hence, we can obtain a Markov partition $\{P_k\}_{k=1}^m$
 associated with $\Phi$
 such that the diameter of each $P_k$ is smaller than
 a given constant $\epsilon>0$ by the same argument
 as the hyperbolic case
 (see {\it e.g.} \cite{Ra} for the hyperbolic case).
By the $C^2$-regularity of $E^s$, we can choose $\{\phi_k\}$
 so that $\pi_y \circ \phi_k^{-1}$ is of class $C^2$,
 where $\pi_y(x,y)=y$.
It implies that $F$ is a piecewise $C^2$ map.
It is easy to check that $F$ is a pre-Markov map
 with the set of discontinuity
 $\left(\bigcup_{k,l=1}^m \del J_{k,l}\right) \- \del I_*$.
By Proposition \ref{prop:stable manifolds 1},
 $F$ is a topologically Markov map.
\end{proof}

For $\sigma=s,u$,
 let $\Per^\sigma_*(\Phi)$ be the set of periodic point $z_*$
 such that $\|N\Phi^{t_*}|_{(E^\sigma/T\Phi)(z_*)}\|=1$,
 where $t_*$ is the period of $z_*$.
The following is an immediate consequence of
 the above lemma and Proposition \ref{prop:Mane}.
\begin{prop}
\label{prop:expansion}
If a $C^2$ $\PA$ flow $\Phi$ is $E^s$-fine,
 $\Per^u_*(\Phi)$ consists of finitely many orbits.
Moreover,
 for any given neighborhood $U$ of $\Per^u_*(\Phi)$,
 there exists $T=T(U)>0$ such that
\begin{equation*}
\sup\left\{\|N\Phi^{-t}|_{(E^u/T\Phi)(z)}\|
 \st t \geq T, z \in M \- U \right\} \leq \frac{1}{2}.
\end{equation*}
\end{prop}
\begin{coro}
\label{cor:Anosov}
If a $C^2$ $\PA$ flow $\Phi$ is $E^s$-fine
 and $\Per^u_*(\Phi)$ is empty,
 then there exists $C>0$ and $\lambda>1$ such that
 $\|N\Phi^{-t}|_{E^u/T\Phi(z)}\|<C\lambda^{-t}$
 for any $z \in M$ and $t >0$.
\end{coro}

%
%
\subsection{Non-expansion property}
\label{sec:non-expansion}
Let $M$ be a closed three-dimensional Riemannian manifold
 and $\Phi$ be a $C^2$ $E^s$-fine $\PA$ flow.
As we saw in Subsection \ref{sec:hyp-like},
 there exists a family $\{V^u(z)\}_{z \in M}$
 of $C^2$ immersed submanifolds of $M$ which are tangent to $E^u$.
For $z \in M$ and $\delta>0$,
 let $B(z,\delta)$ be the closed $\delta$-ball
 centered at $z$
 and $D^u(z,\delta)$ the connected component
 of $V^u(z) \cap B(z,\delta)$ that contains $z$.
By the definition of $V^u(z)$,
 we can see that
\begin{itemize}
 \item $V^u(\Phi^t(z))=\Phi^t(V^u(z))$ for any $z \in M$ and $t \in \RR$,
 \item $D^u(z,\delta)$ is a $C^2$ embedded disk
 which varies continuously with respect to $z$
 if $\delta$ is sufficiently small.
\end{itemize}
We say that $\Phi$ is {\it $u$-bounded}
 if there exists a positive-valued function
 $\bar{\delta}$ on $\RR$ such that
\begin{equation*}
\Phi^{-t}(D^u(z,\bar{\delta}(\epsilon))) \subset D^u(\Phi^{-t}(z),\epsilon) 
\end{equation*}
 for any $z \in M$, $t \geq 0$, and $\epsilon>0$.
Similarly, we can define a $C^2$ disk $D^s(z,\delta)$
 which is tangent to $E^s$ for any $z \in M$
 and any sufficiently small $\delta>0$.
We say that $\Phi$ is {\it $s$-bounded}
 if there exists a positive-valued function
 $\bar{\delta}'$ on $\RR$ such that
\begin{equation*}
\Phi^t(D^s(z,\bar{\delta}'(\epsilon))) \subset D^s(\Phi^t(z),\epsilon) 
\end{equation*}
 for any $z \in M$, $t \geq 0$, and $\epsilon>0$.

If $E^u$ admits a continuous $D\Phi$-invariant splitting
 $E^u=T\Phi \oplus E^{uu}$,
 then we can show $\Phi$ is $u$-bounded
 by Proposition \ref{prop:stable manifolds 1}.
However, a $\PA$ flow is not $u$-bounded in general.
In fact, if there exists a $\Phi$ invariant embedded annulus
 tangent to $E^u$ on which $\Phi^t$ is conjugate to
 the map $(x,y) \mapsto (x+(1+y)t,y)$ on $S^1 \times [0,1]$,
 then $\Phi$ is not $u$-bounded.
 
We say that $\Phi$ admits a local invariant foliation $\cG$
 transverse to a compact $\Phi$-invariant set $\Lambda$
 if $\cG$ is a $C^2$ two dimensional foliation on
 an open neighborhood $U_*$ of $\Lambda$
 which is transverse to the orbits of $\Phi$
 and satisfies $D\Phi^t(T\cG(z))=T\cG(\Phi(z))$
 for any $t \geq 0$ and $z \in \bigcap_{t' \in [0,t]}\Phi^{-t'}(U_*)$.
Remark that a suitable time change of $\Phi$ admits
 a local invariant foliation transverse to $\Per^u_*(\Phi)$
 if $\Per^u_*(\Phi)$ consists of finite number of orbits.

The aim of this subsection is to show 
 the $u$-boundedness of $\Phi$ under a mild assumption.
\begin{prop}
\label{prop:non-expansion}
If a $C^2$ $E^s$-fine $\PA$ flow $\Phi$
 admits a local invariant foliation transverse to $\Per^u_*(\Phi)$,
 then $\Phi$ is $u$-bounded.
\end{prop}

Fix a local invariant foliation $\cG$
 which is transverse to $\Per^u_*(\Phi)$
 on a neighborhood $U_*$ of $\Per^u_*(\Phi)$.
\begin{lemm}
\label{lemma:non-expansion 1}
For any given neighborhood $U$ of $\Per^u_*(\Phi)$,
\begin{equation*}
 \sup\left\{\|D\Phi^{-t}|_{E^u(z)}\| \st t \geq 0, z \in M \- U
 \right\} <\infty.
\end{equation*}
\end{lemm}
\begin{proof}
By taking a finite covering,
 we may assume that $E^u$ is orientable without loss of generality.
By $X$, we denote the vector field generating $\Phi$.
Fix a continuous vector field $Y^u$
 such that $\{X(z),Y^u(z)\}$ spans $E^u(z)$ for any $z \in M$
 and $Y^u(z)$ is tangent to $\cG(z)$ if $z \in U_*$.
We define functions $\eta$ and $\lambda$ on $M \times \RR$ by
\begin{equation*}
D\Phi^{-t}(Y^u(z))=\eta(z,t) X(\Phi^{-t}(z))+\lambda(z,t) Y^u(\Phi^{-t}(z)).
\end{equation*}
Since $D\Phi^{-t}(X(z))=X(\Phi^{-t}(z))$, the following identities hold:
\begin{align*}
\eta(z,t+t') &= \eta(z,t)+\lambda(z,t)\cdot \eta(\Phi^{-t}(z),t'),\\
\lambda(z,t+t') & =\lambda(z,t) \cdot \lambda(\Phi^{-t}(z),t').
\end{align*}
By the local invariance of $\cG$,
 if $\Phi^{-t}(z) \in U_*$ for any $0 < t <t_0$
 then $\eta(z,t_0)=0$.

Without loss of generality, we may assume that
 $U$ is a subset of $U_*$.
By Proposition \ref{prop:expansion},
 there exists $T>0$ such that
 $|\lambda(z,t)|<1/2$ for any $t \geq T$ and $z \in M \- U$.
Take a constant $K_0>1$ such that
 $|\lambda(z,t)|+|\eta(z,t)| \leq K_0-1$ for any $z \in M$
 and any $0 \leq t \leq T$.
It is sufficient to show that
 $|\lambda(z,t)|+|\eta(z,t)| \leq 4 K_0$
 for any $z \in M \- U$ and $t >0$.

We define a function $\tau:M \- U \ra [T,\infty]$ by
\begin{equation*}
\tau(z)=\left\{
\begin{array}{ll}
\inf\{t \geq T \st \Phi^{-t}(z) \not\in U\} & 
 \mbox{ if } z \not\in \bigcap_{t \geq T}\Phi^t(U)\\
\infty & \mbox{otherwise}.
\end{array}
\right.
\end{equation*}
We claim that $|\lambda(z,t)|+|\eta(z,t)| \leq K_0$
 for any $z \in M \- U$ and $0 \leq t \leq \tau(z)$.
It is trivial if $t \leq T$.
Suppose that $T < t \leq \tau(z)$.
Then, $\Phi^{-t'}(z)$ is contained in $U$
 for any $T \leq t'<t$,
 and hence, $\eta(\Phi^{-T}(z),t-T)=0$.
Since $\lambda(z,t) < 1/2$, we have
\begin{align*}
|\lambda(z,t)|+|\eta(z,t)|
 & = |\lambda(z,t)| +\left|\eta(z,T)
   + \lambda(z,T) \cdot \eta(\Phi^{-T}(z),t-T)\right|\\
 & < \frac{1}{2} +|\eta(z,T)| \leq K_0.
\end{align*}
It completes the proof of the claim.

Fix $z \in M \- U$
Take a sequence $(t_i)_{i \geq 0}$ in $[0,\infty]$ such that
 $t_0=0$, $t_{i+1}=t_i+\tau(\Phi^{-t_i}(z))$ for any $i \geq 0$.
We claim 
\begin{equation}
\label{eqn:non-expansion 1}
\lambda(z,t_i) \leq 2^{-i},\hsp
\eta(z,t_i) \leq 2K_0\left(1-2^{-i}\right)
\end{equation}
 for any $i \geq 0$.
The proof is by induction.
The inequalities for $i=0$ is trivial.
Suppose that they hold for $i$.
Since $T \leq t_{i+1}-t_i =\tau(\Phi^{-t_i}(z))$
 and $\tau(\Phi^{-t_i}(z)) \not\in U$,
 we have $|\lambda(\Phi^{-t_i}(z),t_{i+1}-t_i)| \leq 1/2$.
The first claim also implies
$|\eta(\Phi^{-t_i}(z),t_{i+1}-t_i)| \leq K_0$.
By the assumption of induction, 
\begin{align*}
|\lambda(z,t_{i+1})|
 & =|\lambda(z,t_i)|\cdot |\lambda(\Phi^{-t_i}(z),t_{i+1}-t_i)|
 \leq 2^{-i} \cdot 2^{-1} =2^{-(i+1)},\\
 |\eta(z,t_{i+1})| &
 \leq |\eta(z,t_i)|+|\lambda(z,t_i)|
 \cdot |\eta(\Phi^{t_i}(z),t_{i+1}-t_i)|\\
 &\leq 2K_0\left(1-2^{-i}\right)+2^{-i}K_0
 =2K_0\left(1-2^{-(i+1)}\right).
\end{align*}
Therefore, the inequalities (\ref{eqn:non-expansion 1}) hold
 for $i+1$. The claim is proved.

Take $t>0$.
There exists $i \geq 0$ such that $t_i \leq t <t_{i+1}$.
Since $0 \leq t-t_i <\tau(\Phi^{-t_i}(z))$,
 the above claims imply
\begin{align*}
 |\lambda(z,t)| &
 =|\lambda(z,t_i)|\cdot |\lambda(\Phi^{-t_i}(z),t-t_i)|
 \leq 1 \cdot K_0 = K_0,\\
 |\eta(z,t)| &
 \leq |\eta(z,t_i)|+|\lambda(z,t_i)|
 \cdot |\eta(\Phi^{-t_i}(z),t-t_i)|\\
 &\leq 2K_0+1 \cdot K_0 =3K_0.
\end{align*}
\end{proof}

\begin{lemm}
\label{lemma:non-expansion 2} 
There exists a neighborhood $U_1$ of $\Per^u_*(\Phi)$
 and a positive-valued function $\bar{\delta}_1$ on $\RR$ such that
\begin{equation*}
\Phi^{-t}(D^u(z,\bar{\delta}_1(\epsilon)))
 \subset D^u(\Phi^{-t}(z),\epsilon)
\end{equation*}
 for any $\epsilon>0$, $t >0$,
 and $z \in \bigcap_{t' \in [0,t]}\Phi^{t'}(U_1)$.
\end{lemm}
\begin{proof}
For an interval $J \subset \RR$ and a subset $S$ of $M$,
 we denote the set $\{\Phi^t(z) \st z \in S, t \in J \}$
 by $\Phi^J(S)$.

Fix a neighborhood $U_1$ of $\Per^u_*(\Phi)$
 and a constant $0<\epsilon'<\epsilon/2$ such that
 $\cN_{\epsilon'}(U_1) \cap \cN_{\epsilon'}(M \- U_*)=\emptyset$.
Take a family $\{\phi_z\}_{z \in M}$ of local cross-sections
 so that $\Im \phi_z \subset \cG(z)$ for any $z \in U_1$.
Let $\{r_z^t\}$ be the family of returns
 and put $I^u_\delta(z)=D^u(z,\delta) \cap \Im \phi_z$.
Since $\{D^u(z,\Delta)\}$ is a continuous family of $C^2$-disks,
 $\{I^s_\delta(z)\}$ is a continuous family of $C^2$-intervals.
By Proposition \ref{prop:stable manifolds 1},
 there exists $\delta'>0$ such that
 $r_z^{-t}(I^u_{\delta'}(z)) \subset I^u_{\epsilon'}(\Phi^{-t}(z))$
 for any $t \geq 0$ and $z \in M$.
Take an open interval $J \subset \RR$ containing $0$ such that
 $\Phi^J(I_{\epsilon'}^u(z)) \subset D^u(z,2\epsilon')$
 for any $z \in M$.
By the invariance of $\cG$,
 we have $r_z^{-t}(I^u_{\delta'(z)})=\Phi^{-t}(I^u_{\delta'}(z))$
 for any $z \in \bigcap_{t' \in [0,t]}\Phi^{t'}(U_1)$.
We also take $\delta>0$ so that
 $D^u(z,\delta) \subset \Phi^J(I^u_{\delta'}(z))$
 for any $z \in M$.
Then, we have
\begin{align*}
\Phi^{-t}(D^u(z,\delta))
 & \subset \Phi^{-t}(\Phi^J(I^u_{\delta'}(z)))
  =\Phi^J (\Phi^{-t}(I^u_{\delta'}(z)))
  =\Phi^J(r_z^{-t}(I^u_{\delta'}(z)))\\
 & \subset \Phi^J(I^u_{\epsilon'}(\Phi^{-t}(z)))
 \subset D^u(\Phi^{-t}(z), 2\epsilon')
  \subset D^u(\Phi^{-t}(z), \epsilon).
\end{align*}
 for any $t>0$ and $z \in \bigcap_{t' \in [0,t]}\Phi^{t'}(U_1)$.
\end{proof}

Now, we prove Proposition \ref{prop:non-expansion}.
Fix $\epsilon>0$.
Let $U_1$ and $\bar{\delta}_1$ be the neighborhood of $\Per^u_*(\Phi)$
 and the function that are given by Lemma \ref{lemma:non-expansion 2}.
There exists a neighborhood $U_2$ of $\Per^u_*(\Phi)$
 and a constant $0 <\epsilon_2 <\epsilon$ such that
 $\cN_{\epsilon_2}(U_2) \cap \cN_{\epsilon_2}(M \- U_1)=\emptyset$.
By Lemma \ref{lemma:non-expansion 1},
\begin{equation*}
 K=1+\sup\left\{\|D\Phi^{-t}|_{E^u(z)}\| \st t \geq 0, z \in M \- U_2
 \right\}
\end{equation*}
 is finite.
Put $\delta_2=\bar{\delta}_1(\epsilon_2K^{-1})$ and
 $\delta=\min\{\delta_2,\epsilon_2K^{-1}\}$.
It is sufficient to show the inclusion
\begin{equation}
\label{eqn:non-expansion *}
\Phi^{-t}(D^u(z,\delta)) \subset D^u(\Phi^{-t}(z),\epsilon) 
\end{equation}
 for any $z \in M$ and $t \geq 0$.

For $z \not\in U_1$, we have
 $D^u(z,\epsilon_2K^{-1}) \cap U_2=\emptyset$,
 and hence,
\begin{equation}
\label{eqn:non-expansion 3}
\Phi^{-t}(D^u(z,\epsilon_2 K^{-1}))
 \subset D^u(\Phi^{-t}(z),\epsilon_2) \subset D^u(\Phi^{-t}(z),\epsilon)
\end{equation}
 for any $t \geq 0$.
It implies the inclusion (\ref{eqn:non-expansion *}) for $z \not\in U_1$.
If $z \in U_1$,
 put $T=\inf\{t>0 \st \Phi^{-t}(z) \not\in U_1\} \in (0,\infty]$.
For $0 \leq t <T$, we have
\begin{equation*}
\Phi^{-t}(D^u(z,\delta_2)) \subset D^u(\Phi^{-t}(z),\epsilon_2 K^{-1})
 \subset D^u(\Phi^{-t}(z),\epsilon). 
\end{equation*}
It implies the inclusion (\ref{eqn:non-expansion *})
  for the case $T=\infty$ or $0< t <T$.
If $T$ is finite, then
 $\Phi^{-T}(D^u(z,\delta_2)) \subset D^u(\Phi^{-T}(z),\epsilon_2K^{-1})$.
Since $\Phi^{-T}(z) \not\in U_1$,
 we have $\Phi^{-(T+t')}(D^u(z,\delta_2))
 \subset D^s(\Phi^{-(T+t')}(z),\epsilon)$
 for any $t'\geq 0$ by (\ref{eqn:non-expansion 3}).
It implies the inclusion (\ref{eqn:non-expansion *}) for
 the case $z \in U_1$ and $t \geq T$.

%
%
\subsection{Hyperbolicity of periodic orbits}
\label{sec:hyp}
The following proposition is the last piece
 of the proof of Proposition \ref{prop:Anosov}.
\begin{prop}
\label{prop:hyperbolicity}
If a $C^2$ $E^s$-fine $\PA$ flow $\Phi$  is $s$- and $u$-bounded,
 then $\Per^u_*(\Phi)$ is empty.
\end{prop}
\begin{proof}
[Proof of Proposition \ref{prop:Anosov}]
Let $\Phi$ be a topologically transitive
 $\PA$ flow with a $C^2$ $\PA$ splitting.
Topological transitivity implies that $\Omega_*$ is empty.
By Proposition \ref{prop:regular point},
 all periodic points are $s$- and $u$-regular.
In particular, $\Phi$ and $\Phi^{-1}$ are $E^s$-fine.

By Proposition \ref{prop:expansion} for $\Phi$ and $\Phi^{-1}$,
 we see that $\Per_*(\Phi)$ consists of finite number of orbits.
Take a time-change $\Phi_1$ of $\Phi$
 which admits a local invariant foliation transverse to $\Per_*(\Phi)$.
Remark that both $\Phi_1$ and $\Phi_1^{-1}$ are and $E^s$-fine.
We apply Propositions \ref{prop:non-expansion}
 and \ref{prop:hyperbolicity} to $\Phi_1$ and $\Phi_1^{-1}$.
The former implies $\Phi_1$ is $s$- and $u$-bounded.
The latter for $\Phi_1$ implies that $\Per_*^u(\Phi_1)$ is empty
 and the same for $\Phi_1^{-1}$ implies that $\Per_*^s(\Phi_1)$ is empty.
Hence, $\Phi_1$ is an Anosov flow by and Corollary \ref{cor:Anosov}.
Since $\Phi$ is a time-change of $\Phi_1$, also $\Phi$ is.
\end{proof}

The rest of the subsection is devoted to
 the proof of Proposition \ref{prop:hyperbolicity}.
Fix a $C^2$ $E^s$-fine $\PA$ flow $\Phi$
 on a closed three-dimensional manifold $M$
 which is $s$- and $u$- bounded.
Let $TM=E^u +E^s$ be a $\PA$ splitting of $\Phi$
 such that $E^s$ generates a $C^2$ foliation $\cF^s$.
Remark that $D^s(z,\delta) \subset \cF^s(z)$
 for any $z \in M$ and $\delta>0$.

Take a family $\{\psi_z\}_{z \in M}$ of $C^2$ embeddings
 from $[-1,1]^3$ to $M$ such that $\psi_z(0,0,0)=z$,
 $\psi_z([-1,1]^2 \times y) \subset \cF^s(\psi_z(0,0,y))$
 for any $y \in [-1,1]$,
 and the map $(z,w) \mapsto \psi_z(w)$
 from $M \times [-1,1]^3$ to $M$ is of class $C^2$.
By $B(z,\delta)$, we denote the closed ball of radius $\delta$
 which is centered at $z$.
There exists $\epsilon_0>0$ such that
$B(z,8\epsilon_0) \subset \Im \psi_z$ for any $z \in M$.
Since $\Phi$ is $s$- and $u$-bounded,
 we can take $\delta_0>0$ so that
 $\Phi^t(D^s(z,\delta_0)) \subset D^s(\Phi^t(z),\epsilon_0)$
 and $\Phi^{-t}(D^u(z,\delta_0)) \subset D^u(\Phi^{-t}(z),\epsilon_0)$
 for any $z \in M$ and $t \geq 0$.

Suppose that $\Per^u_*(\Phi)$ is non-empty
 and contains a point $p$.
There exists a continuous injective map $H:(-1,1)^2 \ra M$ such that
\begin{enumerate}
\item $H(0,0)=p$,
\item $\Im H \subset \Im \psi_p$,
\item $H(x,\cdot)$ is of class $C^2$ and
$H(x \times (-1,1)) \subset D^u(H(x,0),\delta_0)$
 for any $x \in (-1,1)$, and
\item $H((-1,1) \times y) \subset D^s(H(0,y),\delta_0)$
 for any $y \in (-1,1)$.
\end{enumerate}
We put $V=\bigcup_{x \in (-1,1)}D^u(H(x,0),\delta_0)$.

By Proposition \ref{prop:expansion},
 $\Per^u_*(\Phi)$ consists of finitely many orbits.
Since $\Phi$ is topologically transitive,
 the union of periodic orbits is a dense subset of $M$.
Hence, $H((-1,1)^2)$ contains a hyperbolic periodic point $q$.
Put $(x_*,y_*)=H^{-1}(q)$.

For $x \in [0,x_*]$ and $t \geq 0$,
 we put $J^u(x,t)=\Phi^{-t}(H(x \times [0,y_*]))$.
Let $V(x,t)$ be the arcwise connected component
 of $\Phi^{-t}(V) \cap B(\Phi^{-t}(H(x,0)),3\epsilon_0)$
 which contains $\Phi^{-t}(H(x,0))$.
Since
\begin{equation*}
 J^u(x,t) \subset \Phi^{-t}(D^u(H(x,0),\delta_0)
 \subset D^u(\Phi^{-t}(H(x,0)),\epsilon_0),
\end{equation*}
 we have $J^u(x,t) \subset V(x,t)$.

Let $\pi_y:\RR^3 \to \RR$ be the map defined by $\pi_y(w,x,y)=y$.
Put $I(x,t)=\pi_y \circ \psi_{\Phi^{-t}(H(x,0))}^{-1}(J^u(x,t))$.
We define a map $h_{x,t}:[0,y_*] \ra I(x,t)$ by
\begin{equation*}
h_{x,t}(y)=\pi_y \circ \psi_{\Phi^{-t}(H(x,0))}^{-1}(\Phi^{-t}(H(x,y))). 
\end{equation*}
Remark that it is a $C^2$ diffeomorphism.

The map $h_{x_*,t} \circ h_{0,0}^{-1}$ can be decomposed in two ways;
\begin{align}
\label{eqn:distortion}
h_{x_*,t} \circ h_{0,0}^{-1} 
 & = \left(h_{x_*,t} \circ h_{x_*,0}^{-1}\right)
  \circ \left(h_{x_*,0} \circ h_{0,0}^{-1}\right) \\
 & = \left(h_{x_*,t} \circ h_{0,t}^{-1}\right)
  \circ \left(h_{0,t} \circ h_{0,0}^{-1}\right).\nonumber
\end{align}
We estimate the distortion of each decomposition.
It will lead us to a contradiction.

\begin{lemm}
\label{lemma:distortion 1}
$\{\dist(h_{x_*,t} \circ h_{x_*,0}^{-1},I(x_*0,0)) \st t \geq 0\}$ is bounded.
\end{lemm}
\begin{proof}
Let $T$ be the period of $q$
 and put $h=h_{x_*,T} \circ h_{x_*,0}^{-1}$.
Then, $I(x_*,T) \subset I(x_*,0)$
 and the map $h$ is $C^2$ conjugate to
 the local return map of $\Phi^{-1}$ on $J^u(x_*,0)$,
Since $q$ is a hyperbolic periodic point,
 there exist $C>0$ and $\lambda \in (0,1)$ such that
 $|h^n(I(x_*,T))|<C\lambda^n$ for any $n \geq 0$.
Take $K>0$ so that $|D(\log |Dh|)(y)|<K$
 for any $y \in I_{x_*,0}$.
Then, we have
\begin{align*}
\dist(h_{x_*,nT} \circ h_{x_*,0}^{-1},I(x_*,0))
 & = \dist(h^n,I(x_*,0)) \\
 & = \sum_{m=0}^{n-1} \dist(h, h^m(I(x_*,0)))\\
 & \leq \sum_{m=0}^{n-1} K\cdot C\lambda^n
 < KC(1-\lambda)^{-1}.
\end{align*}
Since $\dist(h_{x_*,t} \circ h_{x_*,0}^{-1},(I(x_*,0)))$
 is continuous with respect to $t$, it is bounded on $[0,T]$.
Hence, the lemma follows from the formula (\ref{eqn:distortion 2}).
\end{proof}

We will show that the distortion of
 the last term of  (\ref{eqn:distortion}) is unbounded.
Once it is shown, it contradicts Lemma \ref{lemma:distortion 1},
 and hence, $\Per_*^u(\Per(\Phi))$ is empty.

\begin{lemm}
\label{lemma:distortion 2}
$\{\dist(h_{0,t} \circ h_{0,0}^{-1},I(0,0)) \st t \geq 0\}$ is unbounded.
\end{lemm}
\begin{proof}
Let $T$ be the period of $p=h_{0,0}^{-1}(0)$.
Put $h=h_{0,T} \circ h_{0,0}^{-1}$.
Since $h$ is $C^2$ conjugate to the return map of $\Phi^{-1}$ on $J^u(0,0)$,
 we have $I(0,T)=h(I(0,0)) \subset I(0,0)$,
 $\bigcap_{n \geq 1}h^n(I(0,0))=\{0\}$, $h(0)=0$,
 and $|Dh(0)|=1$.
By the formula (\ref{eqn:distortion 3}),
\begin{align*}
 \dist(h_{0,nT} \circ h_{0,0}^{-1},I(0,0))
 & =\dist(h^n,I(0,0))\\
 & \geq \left|\log|Dh^n(0)|-\log|h^n(I(0,0))|+\log|I(0,0)|\right|\\
 & =\left|-\log|h^n(I(0,0))|+\log|I(0,0)|\right|.
\end{align*}
The last term goes to infinity as $n$ tends to infinity
 since $\lim_{n \ra \infty}|h^n(I(0,0))|=0$.
\end{proof}

To estimate the distortion of $h_{x_*,t} \circ h_{0,t}^{-1}$,
 we need some preparations.
\begin{lemm}
\label{lemma:distortion 3}
If $V(x_1,t) \cap V(x_2,t) \neq \emptyset$
 for $0 \leq x_1<x_2 \leq x_*$ and $t \geq 0$, then
\begin{equation*}
 \bigcup_{x \in [x_1,x_2]}J^u(x,t)
 \subset \Im \psi_{\Phi^{-t}(H(x_1,t))}
\end{equation*}
\end{lemm}
\begin{proof}
Since $V(x_1,t) \cup V(x_2,t)$ is arcwise connected,
 we can take a continuous map
 $L:[0,1] \ra V(x_1,t) \cup V(x_2,t)$ such that
 $L(0)=\Phi^{-t}(H(x_1,0))$ and $L(1)=\Phi^{-t}(H(x_2,0))$.
It defines a continuous map $l:[0,1] \ra (-1,1)$ such that
 $L(\xi) \in \Phi^{-t}(D^u(H(l(\xi)),0),\delta_0)$
 for any $\xi \in [0,1]$.
The diameter of $V(x_1,t) \cup V(x_2,t)$ is not greater than
 $6\epsilon_0$ and
\begin{equation*}
 J^u(l(\xi),t) \cup \{L(\xi)\}
 \subset \Phi^{-t}(D^u(H(l(\xi),0)),\delta_0)
 \subset D^u(\Phi^{-t}(H(l(\xi),0)), \epsilon_0).
\end{equation*}
Hence, $J^u(l(\xi),t)$
 is contained in $B(\Phi^{-t}(H(x_1,0)),8\epsilon_0)$
 for any $\xi \in [0,1]$.
Since $[x_1,x_2] \subset \Im l$
 and $B(\Phi^{-t}(H(x_1,0)),8\epsilon_0)
 \subset \Im \psi_{\Phi^{-t}(H(x_1,0))}$,
 the proof is complete.
\end{proof}

Let $\Vol(\cdot)$ be the volume on $M$
 associated with the fixed Riemannian
 metric of $M$.
\begin{lemm}
\label{lemma:distortion 4} 
There exists a constant $K_*>0$ such that
\begin{equation*}
|I(x,t)| \leq K_* \Vol(V(x,t))
\end{equation*}
 for any $x \in [0,x_*]$ and $t \geq 0$.
\end{lemm}
\begin{proof}
Since $H([0,x_*] \times [0,y_*])$ is contained in $\Int V$,
 there exists $\epsilon_1>0$ such that
 $D^s(z,\epsilon_1) \subset V$ for any $z \in H([0,x_*] \times [0,y_*])$.
Since $\Phi$ is $s$- and $u$-bounded,
 we can take $\delta_1>0$ such that
 $\Phi^t(D^s(z,\delta_1)) \subset D^s(\Phi^t(z),\epsilon_1)$
 and $\Phi^{-t}(D^u(z,\delta_1)) \subset D^u(\Phi^{-t}(z),\epsilon_1)$
 for any $z \in M$.

Put $C(x,t)=\bigcup_{z \in J^u(x,t)}D^s(z,\delta_1)$.
For $z=\Phi^{-t}(H(x,y)) \in J^u(x,t)$,
$\Phi^t(D^s(z,\delta_1))$ is contained in 
 $D^s(\Phi^t(z),\epsilon_1) = D^s(H(x,y),\epsilon_1)$,
 and hence, in $V$.
Since
\begin{equation*}
J^u(x,t) 
 \subset \Phi^{-t}(D^u(H(x,0),\delta_0))
 \subset  D^u(\Phi^{-t}(H(x,0)),\epsilon_0),
\end{equation*}
 we have
\begin{align*}
 C(x,t)
 & \subset \Phi^{-t}(V) \cap B(\Phi^{-t}(H(x,0)),\epsilon_0+\delta_1)
 \subset V(x,t).
\end{align*}

By the $C^2$ smoothness of the map $(z,w,x,y) \mapsto \psi_z(w,x,y)$,
 there exists $K_0>0$ such that
 $K_0^{-1}\|v\| \leq \|D\psi_z^{-1}(v)\| \leq K_0\|v\|$
 for any $z \in M$, $z' \in \Im \psi_z$, and $v \in T_{z'}M$.
Let $\mbox{Leb}_n$ be the Lebesgue measure on $\RR^n$.
Since
\begin{equation*}
\mbox{Leb}_2\left(\psi_{\Phi^{-t}(H(x,0))}^{-1}
  \left(D^s(z,\delta_1)\right)\right)
 \geq \pi \delta_1^2K_0^{-2}
\end{equation*}
 for any $z \in J^u(x,t)$, we have
\begin{align}
 \label{eqn:lemma-distortion 4}
 |I(x,t)| &=|\pi_y \circ \psi_{\Phi^{-t}(H(x,0))}^{-1}(J^u(x,t))| \\
 & \leq \frac{K_0^2}{\pi\cdot \delta_1^2} \cdot
 \mbox{Leb}_3\left(\psi_{\Phi^{-t}(H(x,0))}^{-1}(C(x,t))\right) \nonumber\\
 & \leq \frac{K_0^5}{\pi\cdot \delta_1^2} \cdot
 \Vol\left(C(x,t)\right)
 \leq \frac{K_0^5}{\pi\cdot \delta_1^2} \cdot
 \Vol\left(V(x,t)\right).
 \nonumber
\end{align}
\end{proof}

Now, we estimate the distortion of $h_{x_*,t} \circ h_{0,t}^{-1}$.
\begin{lemm}
\label{lemma:distortion 5}
$\{\dist(h_{x_*,t} \circ h_{0,t}^{-1},I(0,t))\} \st t \geq 0\}$ is bounded.
\end{lemm}
\begin{proof}
For $z,z' \in M$,
 the map $\psi_z^{-1} \circ \psi_{z'}$ can be written as
\begin{equation*}
\psi_z^{-1} \circ \psi_{z'}^{-1}(w,x,y)
 =(f_{z,z'}(w,x,y),g_{z,z'}(y))
\end{equation*}
 by a map $f_{z,z'}$ valued in $\RR^2$ and a function $g_{z,z'}$.
Take $K_1>0$ so that $|D(\log|Dg_{z,z'}|)(y)|\leq K_1$
 for any $z,z' \in M$ and any $y$ in the domain of $g_{z,z'}$.

By Lemma \ref{lemma:distortion 3},
 if $V(x_1,t) \cap V(x_2,t) \neq \emptyset$
 for $0 \leq x_1 < x_2 \leq x_*$, then
 $\pi_y \circ \psi_{\Phi^{-t}(H(x_1,0))}(J^u(x,t))=I(x_1,t)$
 for any $x \in [x_1,x_2]$.
It implies that
$h_{x_2,t} \circ h_{x_1,t}^{-1}=g_{z_2(t),z_1(t)}$,
 where $z_i(t)=\Phi^{-t}(H(x_i,0))$.
Hence, we have
\begin{equation*}
 \dist(h_{x_2,t} \circ h_{x_1,t}^{-1},I(x_1,t)) \leq K_1.
\end{equation*}

Fix $t \geq 0$.
Let $\cS$ be the set of sequences $(x_i)_{i=0}^m$ that satisfy
 $x_0=0$, $x_m=x_*$, and
 $V(x_{i+1},t) \cap V(x_i,t) \neq \emptyset$ for any $i=0,\cdots,m-1$.
It is non-empty by the compactness of $[0,1]$.
Take $(x_i)_{i=0}^m \in \cS$ such that $m$ is minimal in $\cS$.
For any $z \in M$,
 the minimality of $m$ implies that
 the number of $V(x_i,t)$ containing $z$ is at most two.
Let $K_*$ be the constant given by Lemma \ref{lemma:distortion 4}.
Then, we have
\begin{align*}
 \dist(h_{x_*,t} \circ h_{0,t}^{-1},I(0,t))
 & \leq \sum_{i=0}^{m-1}
  \dist(h_{x_{i+1},t} \circ h_{x_i,t}^{-1},I(x_i,t))\\
 & \leq K_1 \sum_{i=0}^{m-1}|I(x_i,t)|\\
 & \leq K_* K_1 \sum_{i=0}^{m-1}\Vol(V(x_i,t))\\
 & \leq 2K_*K_1 \Vol(M).
\end{align*}
Since $K_1$ and $K_*$ does not depend on $t$,
 the lemma is proved.
\end{proof}

Since the second component of the middle term of  (\ref{eqn:distortion})
 does not depend on $t$,
 Lemma \ref{lemma:distortion 1} implies
 that the distortion of the middle term is bounded with respect to $t$.
It contradicts
 Lemmas \ref{lemma:distortion 2} and \ref{lemma:distortion 5},
 which imply that the distortion of
 the last term of (\ref{eqn:distortion})
 is unbounded with respect to $t$.
Therefore, $\Per_*^u(\Phi)$ is empty.
Now, the proof of Proposition \ref{prop:hyperbolicity} is finished.

\section{Foliations with tangentially contracting flows}
\label{sec:GA}

%
%
\label{sec:TC}
In this section, we prove Theorem \ref{thm:TC}.

Let $\cF$ be a $C^2$ foliation
 on a closed three-dimensional manifold $M$.
Suppose that $\cF$ admits a $C^2$ tangentially contracting flow $\Phi$.
Let $C>0$ and $\lambda>1$ be constants such that
 $\|N\Phi^t|_{(T\cF/T\Phi)(z)}\| \leq C\lambda^{-t}$
 for any $z \in M$ and $t \geq 0$.

\begin{lemm}
\label{lemma:TC is PA}
There exists a continuous subbundle $E^u$ of $TM$
 such that $\Phi$ is a $\PA$ flow with a $\PA$ splitting
 $TM=T\cF + E^u$.
\end{lemm}
\begin{proof}
The proof is almost identical to Lemme IV.1.1 in \cite{Gh2}.

The differential of the flow $\Phi$
 induces a flow $N_\cF\Phi$ on $TM/T\cF$.
Let $S_*$ be the set of points $z \in M$ that satisfies
\begin{equation}
\label{eqn:TC PA}
 \limsup_{t \ra \infty}\frac{1}{t}
 \log \|N_\cF\Phi^t_z\|\leq -\frac{\log\lambda}{3}.
\end{equation}
We will show that $S_*$ must be empty.
Once it is shown, we have
\begin{equation*}
 \liminf_{t \ra \infty}
 \frac{\|N\Phi^t_{(T\cF/T\Phi)(z)}\|}{\|N_\cF\Phi^t_z\|} =0
\end{equation*}
 for any $z \in M$.
By a standard argument (see {\it e.g. } Proposition 2.3 of \cite{As3}),
 we can show that there exists a continuous subbundle $E^u$ of $TM$
 such that $TM=T\cF + E^u$ is a $\PA$ splitting for $\Phi$.

Suppose that $S_*$ is not empty.
Take a point $z_0$ in $S_*$.
First, we claim that $\Phi$ admits an attracting periodic point.
As an accumulation point of the uniform measures
 on $\{\Phi^t(z_0) \st t \in [0,T]\}$ with $T \ra \infty$,
 we obtain a $\Phi$-invariant Borel probability measure $m_*$
 such that
\begin{displaymath}
\int_M \left(\left.\frac{d}{dt}
 \log\|N_\cF\Phi^t_z\|\right|_{t=0}\right) dm_*(z)
 \leq -\frac{\log\lambda}{3}.
\end{displaymath}
It implies that there exists a Borel subset $U_*$ of $M$
 such that $m_*(U_*)>0$
 and all Lyapunov exponents of $\Phi$ are negative on $U_*$.
By Pesin theory, the $\omega$-limit set of any point of $U_*$
 is an attracting periodic point.
Therefore, the claim is proved.

Suppose that $z_a$ is an attracting periodic point of $\Phi$.
Since $\Phi$ is tangentially contracting,
 there exists a compact embedded annulus $A$ in $\cF(z_a)$
 such that $\Phi^t(A) \subset \Int A$ for any $t>0$,
 $\bigcap_{t \geq 0}\Phi^t(A)=\cO(z_a)$, and
 $\bigcup_{t \geq 0}\Phi^{-t}(A)=\cF(z_a)$.
In particular, the leaf $\cF^s(z_a)$ is diffeomorphic to
 $S^1 \times \RR$.
Since $z_a$ is attracting,
 we can take a compact neighborhood $U$ of $\cO(z_a)$ in $M$
 such that $\del A \cap U=\emptyset$ and $\Phi^t(U) \subset \Int U$
 for any $t > 0$.
By the choice of $A$, we have $U \cap \cF(z_a)=U \cap A$.
It implies that $\cF(z_a)$ is a proper leaf.
By Lemma \ref{lemma:semi-proper}, $\cF^s(z_a)$ has trivial holonomy.
However, it contradicts that $\cO(z_a)$ is an attracting periodic orbit.
\end{proof}

\begin{lemm}
\label{lemna:TC transitive}
The $\PA$ flow $\Phi$ is $E^s$-fine.
\end{lemm}
\begin{proof}
By the strong stable manifold theorem,
 each leaf of $\cF$ is diffeomorphic to $\RR^2$ or $S^1 \times \RR$.
In particular,  $\cF$ has no closed leaves.
By Duminy's theorem, there exists no exceptional minimal set of $\cF$.
Hence, each leaf of $\cF$ is dense in $M$.
By the same argument as the above lemma,
 if $z_0$ is a $u$-irregular periodic point,
 then $\cF(z_0)$ is semi-proper.
However, it contradicts that each leaf of $\cF$ is dense in $M$.
Therefore, any periodic point is $u$-regular.
Since $\Phi$ is tangentially contracting,
 any periodic point is $s$-regular.

By Proposition \ref{prop:reduced dichotomy},
 either $\Omega_*$ is empty or $M=W^u(\Omega_*^s) \cap \Omega_*^u$.
Since $\cF$ has no closed leaves, $\Omega_*^s$ is empty.
It implies that the latter case can not occur.
Therefore, $\Omega_*$ is empty.
\end{proof}

\begin{prop}
The flow $\Phi$ is an Anosov flow.
\end{prop}
\begin{proof}
Let $TM=T\cF +E^u$ be a $\PA$ splitting for $\Phi$.
Since $\Phi$ is tangentially contracting with respect to $\cF$,
 we have $\Per^s_*(\Phi)=\emptyset$.
In particular, $\Per_*(\Phi)=\Per_*^u(\Phi)$.

By Proposition \ref{prop:expansion},
 $\Per^u_*(\Phi)$ consists of finitely
 many non-hyperbolic periodic orbits if it is not empty.
Any time-change of $\Phi$ preserves each leaf of $\cF$ and
 is tangentially contracting with respect to $\cF$.
Hence, we may assume that $\Phi$ admits
 a local invariant foliation transverse to $\Per_*(\Phi)=\Per_*^u(\Phi)$
 (see the beginning of Subsection \ref{sec:hyp}).
By Proposition \ref{prop:non-expansion}, $\Phi$ is $u$-bounded.

By the same argument as the hyperbolic case in \cite{Do},
 there exists a continuous $D\Phi$-invariant splitting
 $T\cF=T\Phi \oplus E^{ss}$ and 
 constants $C>0$ and $\lambda>0$ such that
 $\|D\Phi^t|_{E^{ss}(z)}\| \leq C\lambda^{-t}$
 for any $z \in M$ and $t \geq 0$.
It implies that
 $\{\|D\Phi^t|_{T\cF(z)}\| \st z \in M, t \geq 0\}$ is bounded.
Hence, $\Phi$ is $s$-bounded.

Since $\Phi$ is a $C^2$ $E^s$-fine $\PA$ flow,
 Proposition \ref{prop:hyperbolicity} implies
 $\Per_*(\Phi)=\Per_*^u(\Phi)=\emptyset$.
By Corollary \ref{cor:Anosov}, $\Phi$ is an Anosov flow.
\end{proof}

Now, we prove Theorem \ref{thm:TC}.
By Th\'eor\`eme 4.1 of \cite{Gh},
 $\cF$ admits a $C^r$ transverse projective structure.
By Th\'eor\`eme 5.1 of \cite{Ba}
 ({\it c.f.} Th\'eor\`eme 4.7 of \cite{Gh})
 $\Phi$ is topologically equivalent to an algebraic Anosov flow.
Therefore, $\cF$ is homeomorphic to
 the weak stable foliation of an algebraic Anosov flow.
Such a $C^r$ foliation with $r \geq 2$ is classified completely
 by Ghys \cite{Gh} and Ghys and Sergiescu \cite{GS}.
Their result implies that $\cF$ is $C^r$-diffeomorphic to
 the weak stable foliation of an algebraic Anosov flow.

%
%


\begin{thebibliography}{10}
\bibitem{AR}
A.~Arroyo and F.~Rodriguez Hertz,
 Homoclinic bifurcations and uniform hyperbolicity
 for three-dimensional flows.
 {\it Ann. Inst. H. Poincar\'e Anal. Non Lin\'eaire} {\bf 20} (2003), 805--841.

\bibitem{As0}
M.~Asaoka,
 Non-homogeneous locally freeactions of the affine group.
 preprint. arXiv:math.0702833. 

\bibitem{As}
M.~Asaoka,
 Classification of regular and non-degenerate
 projectively Anosov diffeomorphisms on three dimensional manifolds.
 {it. J. Math. Kyoto Univ.} {\bf 46} (2006), no. 2, 349--356.

\bibitem{As2}
M.~Asaoka,
A classification of three dimensional regular projectively Anosov flows.
{\it Proc. Japan Acad. Ser. A Math. Sci.} {\bf 80} (2004),
 no. 10, 194--197 (2005).

\bibitem{As3}
M.~Asaoka,
 Codimension-one foliations with a transversely contracting flow.
{\it Foliations 2005}, 21--36,
 {\it World Sci. Publ., Hackensack, NJ}, 2006.

\bibitem{Ba}
T.~Barbot,
Caract\'erisation des flots d'Anosov en dimension 3
 par leurs feuilletages faibles.
{\it Ergodic Theory Dynam. Systems} {\bf 15} (1995), no. 2, 247--270.

\bibitem{BDV}
C.~Bonatti, L.~J.~D\'iaz, and M.~Viana,
Dynamics beyond uniform hyperbolicity.
A global geometric and probabilistic perspective.
 Encyclopaedia of Mathematical Sciences, {\bf 102}.
 Mathematical Physics, III.
{\it Springer-Verlag, Berlin}, 2005.

\bibitem{CC}
J.~Cantwell and L.~Conlon,
The theory of levels.
 {\it Index theory of elliptic operators, foliations,
 and operator algebras (New Orleans, LA/Indianapolis, IN, 1986)},
 1--10, Contemp. Math., 70, {\it AMS, Providence, RI}, 1988. 

\bibitem{CC1}
J.~Cantwell and L.~Conlon,
Reeb stability for noncompact leaves in foliated 3-manifolds.
 {\it Proc.~Amer.Math.Soc.}, {\bf 33}, 1981, no.~2, 408--410. 

\bibitem{CC2}
J.~Cantwell and L.~Conlon,
 Endsets of exceptional leaves; a theorem of G.Duminy.
{\it Foliations: geometry and dynamics (Warsaw, 2000)}, 225--261, 
{\it World Sci. Publishing, River Edge, NJ}, 2002. 

\bibitem{Do}
C.~I.~Doering,
Persistently transitive vector fields on three-dimensional manifolds.
{\it Dynamical systems and bifurcation theory (Rio de Janeiro, 1985)},
 59--89,
Pitman Res. Notes Math. Ser., {\bf 160},
{\it Longman Sci. Tech., Harlow}, 1987. 

\bibitem{ET}
  Y.~Eliashberg and W.~Thurston,
  Confoliations. University Lecture Series, {\bf 13}. 
  {\it Amer.~ Math.~ Soc., Providence, RI}, 1998.

\bibitem{Gh2}
 E.~Ghys,
 Actions localement libres du groupe affine.
 {\it Invent.~Math.} {\bf 82} (1985), no.3, 479--526.

\bibitem{Gh}
 E.~Ghys, Rigidit\'e diff\'erentiable des groupes fuchiens.
 {\it Inst. Hautes \'Etudes Sci. Publ. Math}.
  {\bf 78} (1993), 163--185.

\bibitem{GS}
E.~Ghys and V.~Sergiescu,
 Stabilit\'e et conjugaison diff\'erentiable pour certains feuilletages.
 {\it Topology} {\bf 19} (1980), no. 2, 179--197.

\bibitem{Ho}
A.~J.~Homburg,
Piecewise smooth interval maps with non-vanishing derivative.
 {\it Ergodic Theory Dynam. Systems} {\bf 20} (2000), no.~3,
 749--773.


\bibitem{Ma}
R.~Ma\~n\'e,
Hyperbolicity, sinks and measure in one-dimensional dynamics.
{\it Comm. Math. Phys.} {\bf 100} (1985), no. 4, 495--524.
 and Erratum. {\it Comm. Math. Phys.} {\bf 112} (1987), no. 4,
 721--724.

\bibitem{Mi}
 Y.~Mitsumatsu,
  Anosov flows and non-stein symplectic manifolds.
  {\it Ann. Inst. Fourier} {\bf{45}} (1995), no.~5,
 1407--1421.

\bibitem{Mi2}
Y.~Mitsumatsu,
  Foliations and contact structures on 3-manifolds.
{\it Foliations: geometry and dynamics (Warsaw, 2000)}, 75--125, 
{\it World Sci. Publishing, River Edge, NJ}, 2002. 

\bibitem{MS}
W.~de Melo and S.~van Strien,
 One-dimensional dynamics.
 Ergebnisse der Mathematik und ihrer Grenzgebiete (3), 25.
 {\it Springer-Verlag, Berlin}, 1993. 

\bibitem{MR}
R.Moussu and R.Roussarie,
Relations de conjugaison et de cobordisme entre certains feuilletages.
{\it Inst. Hautes \'Etudes Sci. Publ. Math.} {\bf 43} (1974), 142--168.

\bibitem{NP}
S.~Newhouse and J.~Palis,
 Hyperbolic nonwandering sets on two-dimensional manifolds.
{\it Dynamical systems (Proc. Sympos., Univ. Bahia, Salvador, 1971)},
 pp. 293--301. Academic Press, New York, 1973. 

\bibitem{No}
 T.~Noda,
 Projectively Anosov flows with differentiable (un)stable foliations.
 {\it  Ann.~Inst.~Fourier} {\bf{50}} (2000), no~.5, 1617--1647.

 \bibitem{No2}
T.~Noda,
 Regular projectively Anosov flows with compact leaves.
 {\it Ann. Inst. Fourier (Grenoble)} {\bf 54} (2004), no. 2, 353--363.

\bibitem{NT}
T.~Noda and T.~Tsuboi,
 Regular projectively Anosov flows without compact leaves.
{\it Foliations: geometry and dynamics (Warsaw, 2000)},
 403--419,
 {\it World Sci. Publishing, River Edge, NJ}, 2002.


\bibitem{PT}
J.~Palis and F.~Takens,
Hyperbolicity and sensitive chaotic dynamics at homoclinic bifurcations.
Cambridge Studies in Advanced Mathematics, {\bf 35}.
{\it Cambridge University Press, Cambridge}, 1993. 

\bibitem{Ra}
M.~Ratner, 
Markov decomposition for an U-flow on a three-dimensional manifold.
 (Russian) {\it. Mat. Zametki} {\bf 6}, 1969, 693--704. 
 (Translation to English) {\it Math. Notes} {\bf 6} 1969, 880--886.

\bibitem{Sh}
M.~Shub,
 Global stability of dynamical systems.
 {\it Springer-Verlag, New York}, 1987.

\bibitem{Ts}
T.~Tsuboi,
 Regular projectively Anosov flows on the Seifert fibered 3-manifolds,
 {\it J. Math. Soc. Japan.}. {\bf 56} (2004), no. 4, 1233--1253.
\end{thebibliography}
\end{document}